\documentclass{imsart}

\RequirePackage{amsthm,amsmath,amsfonts,amssymb}
\RequirePackage[authoryear]{natbib}%% uncomment this for author-year citations
\RequirePackage[colorlinks,citecolor=blue,urlcolor=blue]{hyperref}%% uncomment this for coloring bibliography citations and linked URLs
\RequirePackage{graphicx}%% uncomment this for including figures
\RequirePackage[capitalize]{cleveref}
\allowdisplaybreaks

\startlocaldefs
\theoremstyle{plain}
\newtheorem{theorem}{Theorem}
\newtheorem{lemma}{Lemma}
\theoremstyle{remark}
\newtheorem{remark}{Remark}

\DeclareMathOperator{\ess}{ess}
\DeclareMathOperator{\cov}{cov}
\DeclareMathOperator{\vc}{vc}
\DeclareMathOperator{\HC}{HC}
\DeclareMathOperator{\BB}{BB}
\DeclareMathOperator{\MT}{MT}
\endlocaldefs

\begin{document}

\begin{frontmatter}
\title{Asymptotics of higher criticism via Gaussian approximation}
%\title{A sample article title with some additional note\thanksref{T1}}
\runtitle{Gaussian approximation for higher criticism}
%\thankstext{T1}{A sample of additional note to the title.}

\begin{aug}
%%%%%%%%%%%%%%%%%%%%%%%%%%%%%%%%%%%%%%%%%%%%%%%
%% Only one address is permitted per author. %%
%% Only division, organization and e-mail is %%
%% included in the address.                  %%
%% Additional information can be included in %%
\author[A]{\fnms{Jingkun}~\snm{Qiu}\ead[label=e1]{jkqiu@stu.pku.edu.cn}}%,
%\author[B]{\fnms{???}~\snm{???}\ead[label=e2]{???@???}}
%\and
%\author[B]{\fnms{???}~\snm{???}\ead[label=e3]{???@???}}
%%%%%%%%%%%%%%%%%%%%%%%%%%%%%%%%%%%%%%%%%%%%%%
%% Addresses                                %%
%%%%%%%%%%%%%%%%%%%%%%%%%%%%%%%%%%%%%%%%%%%%%%
\address[A]{Guanghua School of Management, Peking University\printead[presep={,\ }]{e1}}

%\address[B]{???\printead[presep={,\ }]{e2,e3}}
\end{aug}

\begin{abstract}
Higher criticism is a large-scale testing procedure that can attain the optimal detection boundary for sparse and faint signals.
However, there has been a lack of knowledge in most existing works about its asymptotic distribution for more realistic settings other than the independent Gaussian assumption while maintaining the power performance as much as possible.
In this paper, we develop a unified framework to analyze the asymptotic distributions of the higher criticism statistic and the more general multi-level thresholding statistic when the individual test statistics are dependent $t$-statistics under a finite ($2+\delta$)-th moment condition, $0<\delta\leq1$.
The key idea is to approximate the global test statistic by the supremum of an empirical process indexed by a normalized class of indicator or thresholding functions, respectively.
A new Gaussian approximation theorem for suprema of empirical processes with dependent observations is established to derive the explicit asymptotic distributions.
\end{abstract}

%\begin{keyword}[class=MSC]
%\kwd[Primary ]{???}
%\kwd{???}
%\kwd[; secondary ]{???}
%\end{keyword}


\end{frontmatter}
%%%%%%%%%%%%%%%%%%%%%%%%%%%%%%%%%%%%%%%%%%%%%%
%% Please use \tableofcontents for articles %%
%% with 50 pages and more                   %%
%%%%%%%%%%%%%%%%%%%%%%%%%%%%%%%%%%%%%%%%%%%%%%
%\tableofcontents

%%%%%%%%%%%%%%%%%%%%%%%%%%%%%%%%%%%%%%%%%%%%%%
%%%% Main text entry area:
%\tableofcontents
\section{Introduction}
\label{sec:intro}
Multiple or large-scale testing problems play a central role in modern statistical analysis, for instance differential gene expression analysis, large language model watermark detection, and climate change detection and attribution analysis.
Each aforementioned problem can be decomposed into a set of individual testing problems, which can produce a sequence of individual $p$-values.
Typically, the goal of the detection problem is to test a global null hypothesis in which the individual null hypotheses are all true and thus the individual $p$-values should be combined in an informative way.
Common practices of combining $p$-values include Bonferroni correction \citep{Bon36}, stagewise procedures controlling familywise error rates \citep{Hol79-MR538597, Hoc88-MR995126, Hom88}, and false discovery rate control procedures \citep{BH95-MR1325392, BY01-MR1869245}, among others.
\cite{DJ04-MR2065195} developed \cite{Tuk76}'s approach to propose the higher criticism procedure that can attain \cite{Ing97-MR1456646}'s optimal detection boundary.
We refer to \cite{DJ15-MR3317751} for a comprehensive exposition of this method.

Let $Y_{1},\dots,Y_{n}$ be i.i.d. $\mathbb{R}^{p}$-valued random vectors with mean $\mathbb{E}Y_{i}=\mu$.
To test the global hypotheses $H_{0}:\mu=0$ versus $H_{1}:\mu\neq0$, we construct test statistic $T_{j}=T(Y_{1j},\dots,Y_{nj})$ for each dimension $j=1,\dots,p$ and obtain the corresponding individual $p$-values denoted by $p_{j}$.
The higher criticism statistic is given by
\begin{equation}
\label{eq:HC-significance-level}
\HC_{p}^{*}=\sup_{\alpha\in[\alpha_{1},\alpha_{2}]}\frac{1}{(p\alpha(1-\alpha))^{1/2}}\sum_{j=1}^{p}(1_{\{p_{j}\leq\alpha\}}-\alpha),
\end{equation}
where $\alpha$ is a candidate significance level and $[\alpha_{1},\alpha_{2}]\subset(0,1)$ is the searching range of candidate significance levels.
When $[\alpha_{1},\alpha_{2}]=[1/p,1/2]$ and $T_{j}$'s are independent $N(0,1)$ variables, it follows from \cite{Jae79-MR0515687} that as $p\to\infty$,
\begin{equation}
\label{eq:HC-significance-level-asymptotics}
\mathbb{P}\Big(\HC_{p}^{*}\leq\frac{2\log\log p+(1/2)\log\log\log p+x}{(2\log\log p)^{1/2}}\Big)\to\exp\Big(-\frac{1}{2\pi^{1/2}}e^{-x}\Big).
\end{equation}
Then the higher criticism test can be formulated.

A key step of proving \eqref{eq:HC-significance-level-asymptotics} is to show that, on a rich enough probability space, there is a sequence of random variables $\BB_{p}^{*}=_{d}\sup\{B_{\alpha}^{0}/(\alpha(1-\alpha))^{1/2}:\alpha\in[\alpha_{1},\alpha_{2}]\}$ such that
\begin{equation}
\label{eq:HC-significance-level-GA}
|\HC_{p}^{*}-\BB_{p}^{*}|=o_{\mathbb{P}}((2\log\log p)^{-1/2}),
\end{equation}
where $B^{0}$ is a Brownian bridge process on $[0,1]$ with $\mathbb{E}B_{t}^{0}=0$ and $\mathbb{E}B_{s}^{0}B_{t}^{0}=s(1-t)$ for every $0\leq s\leq t\leq 1$.
In the proof of \cite{Jae79-MR0515687}, the coupling type Gaussian approximation error bound \eqref{eq:HC-significance-level-GA} was ensured by a strong approximation result in \cite{KMT75-MR375412}.
In a seminal paper, \cite{CCK14-MR3262461} generalized \eqref{eq:HC-significance-level-GA} for suprema of empirical processes indexed by classes of functions.
The results were established for i.i.d. observations, which corresponds to i.i.d. individual $p$-values in the higher criticism and is thus too restrictive for this purpose.

In this paper, we shall extend \cite{CCK14-MR3262461}'s Gaussian approximation result for dependent observations.
As an application, we shall further derive the asymptotic distribution of the higher criticism statistic for dependent Gaussian and non-Gaussian test statistics $T_{j}$, especially Student's $t$-statistics.
To the best of our knowledge, the asymptotic analysis of the higher criticism statistic under dependence and/or non-Gaussianity is largely an open question in the high dimensional statistical inference literature.
\cite{HJ10-MR2662357} proposed an innovated higher criticism statistic via data transformation when $T_{j}$'s are dependent Gaussian variables without providing the asymptotic distribution.
Let $p=n^{1/\theta}$ with $\theta\in(0,1)$.
\cite{DHJ11-MR2815777} considered the bootstrap higher criticism $t$-statistic without the explicit asymptotic distribution when $[\alpha_{1},\alpha_{2}]=[1/p,(\log p)/p^{1-\theta}]$ and $T_{j}$'s are independent $t$-statistics.
In a pioneering work, \cite{ZCX13-MR3161449} characterized the asymptotic distributions of the higher criticism statistic and the more general multi-level thresholding statistic when $[\alpha_{1},\alpha_{2}]=[L_{p}/p^{1-\eta},L_{p}/p^{1-\theta}]$ and $T_{j}$'s are dependent $z$-statistics for an arbitrarily small constant $\eta>0$ and some polylogarithmic terms $L_{p}$.
However, a careful analysis in \cref{sec:DB} reveals that the exclusion of the range $\alpha\in[L_{p}/p^{1-\theta},1/2]$ causes a slight power loss for moderately sparse signals while the exclusion of the range $\alpha\in[1/p,L_{p}/p^{1-\eta}]$ results in a severe power loss for the highly sparse signal regime.
\cite{LX19-MR3941262} derived a Gaussian approximation result for the higher criticism statistic in a regression setting and only considered the smallest $c\log p$ $p$-values, which can be approximately translated into the range of $[\alpha_{1},\alpha_{2}]=[1/p,c(\log p)/p]$ under $H_{0}$ and is thus too restrictive, as the region is negligible for i.i.d. $p$-values in view of Lemma 4 in \cite{Jae79-MR0515687}.

The key idea in this paper is to observe that the higher criticism statistic can be represented as or approximated by the supremum of an empirical process indexed by a class of normalized indicator functions.
Then the explicit asymptotic distribution can be derived by applying the newly developed Gaussian approximation result for dependent observations.
When the individual test statistics $T_{j}$ are dependent Gaussian variables, we establish the first asymptotic distribution of the higher criticism statistic with the range of $[\alpha_{1},\alpha_{2}]=[(\log p)^{c}/p,1/(\log p)^{d}]$, which only differs from the range of $[\alpha_{1},\alpha_{2}]=[1/p,1/2]$ in \eqref{eq:HC-significance-level-asymptotics} up to some polylogarithmic terms and is enough to attain \cite{Ing97-MR1456646}'s optimal detection boundary.
We also consider the case of $T_{j}$'s being dependent $t$-statistics with $[\alpha_{1},\alpha_{2}]=[(\log p)^{c}/p,1/p^{d}]$, where $d\in((1-\theta)\vee\rho_{0}^{2},1)$ and $\rho_{0}$ controls the correlations among different testing problems.
Although the exclusion of $\alpha\in[1/p^{d},1/2]$ is required to address the non-Gaussianity, we do not need to omit the range of $\alpha\in[(\log p)^{c}/p,L_{p}/p^{1-\eta}]$ that is crucial for highly sparse signal detection so that the higher criticism $t$-test can still have a reasonable power performance as indicated by the analysis in \cref{sec:DB}.

The theoretical framework presented in this paper does not limited to the higher criticism statistic or mean testing problems.
To appreciate this, we extend the asymptotic results for the multi-level thresholding statistic in \cite{ZCX13-MR3161449} and \cite{CLZ19-MR3911118}, which can be approximated by the supremum of an empirical process indexed by a class of normalized thresholding functions.
With the help of a recent bivaraite Cram\'er type moderate deviation result for self-normalized sums in \cite{QCS25+}, we can reduce the exponential type moment conditions in the mentioned papers to a finite ($2+\delta$)-th moment condition with $0<\delta\leq1$ while maintaining the power performance for highly sparse signals.
When the underlying scientific problem is not translated into a mean testing problem, similar asymptotic results may still be established as long as the corresponding Cram\'er type moderate deviation results are available, for instance \cite{SZ16-MR3498022} for Studentized $U$-statistics and more general self-normalized processes.

The rest of this paper is organized as follows.
\cref{sec:GA} provides the main Gaussian approximation theorem for empirical processes with dependent observations.
\cref{sec:HC} establishes the asymptotic distribution of the higher criticism statistic for dependent Gaussian and non-Gaussian individual test statistics.
In particular, \cref{sec:HC-t} considers the higher criticism $t$-test under dependence.
\cref{sec:MTT} further extends the asymptotic results for the multi-level thresholding $t$-test.
The proofs of asymptotic results, Gaussian approximation theorems, and preliminary lemmas can be found, respectively, in Sections \ref{pf:thm-GA-HC}--\ref{pf:thm-GA-MT-t}, \ref{pf:lem-CCK14}--\ref{pf:thm-GA-mixing}, and \ref{pf:lem-comparison-process}--\ref{pf:lem-VC-cov-residual} in sequence.
\cref{sec:DB} provides an additional power analysis on the effect of the candidate significance level searching range $[\alpha_{1},\alpha_{2}]$ on the higher criticism test that has been mentioned above.

\section{Gaussian approximation theorems}
\label{sec:GA}
We begin with an introduction of several notations.
For a measurable space $(\mathcal{X},\mathcal{A})$, a probability measure $Q$, and a measurable function $f:\mathcal{X}\to\mathbb{R}$, we denote $Qf=\int fdQ$, $\|f\|_{Q,r}=(Q|f|^{r})^{1/r}$ for $1\leq r<\infty$, and $\|f\|_{Q,\infty}=Q\text{-}\ess\sup f$ for $r=\infty$.
The convention that $1/q=0$ is adopted for $q=\infty$.
Let $\ell^{\infty}(\mathcal{F})$ be the space of all uniformly bounded functions $Q:\mathcal{F}\to\mathbb{R}$ such that $\|Q\|_{\mathcal{F}}:=\sup\{|Qf|:f\in\mathcal{F}\}<\infty$.
A class $\mathcal{F}$ of functions $f:\mathcal{X}\to\mathbb{R}$ is said to be pointwise measurable if $\mathcal{F}$ contains a countable subset $\mathcal{G}$ such that for every $f\in\mathcal{F}$ there is a sequence $g_{m}\in\mathcal{G}$ with $g_{m}(x)\to f(x)$ for every $x\in \mathcal{X}$.
A nonnegative measurable function $F:\mathcal{X}\to\mathbb{R}$ is said to be an envelope function of $\mathcal{F}$ if $|f(x)|\leq F(x)$ for every $x\in \mathcal{X}$ and $f\in\mathcal{F}$.
The Vapnik--\v{C}ervonenkis (VC) dimension $\vc(\mathcal{F})$ of a class $\mathcal{F}$ of functions $f:\mathcal{X}\to\mathbb{R}$ is the largest number $n$ such that we can find points $(x_{1},t_{1}),\dots,(x_{n},t_{n})$ in $\mathcal{X}\times\mathbb{R}$ shattered by the subgraphs $\{(x,t):t<f(x)\}$ of $\mathcal{F}$, namely for every subset $I\subset\{1,\dots,n\}$, there is a function $f\in\mathcal{F}$ with $t_{i}<f(x_{i})$ for $i\in I$ and $t_{i}\geq f(x_{i})$ for $i\notin I$.
Unless explicitly stated, the underlying probability space throughout of this paper is assumed to be rich enough to admit a sequence of independent uniform random variables on $(0,1)$.
Let $\delta_{x}$ denote the Dirac measure at $x\in\mathcal{X}$.

We first state a Gaussian approximation theorem for empirical processes indexed by VC classes of functions with i.i.d. observations, which is modified from Corollary 2.2 in \cite{CCK14-MR3262461}.
Although there are other types of results in the mentioned paper, it turns out that the VC type Gaussian approximation results are the most useful in practice, especially the higher criticism as will be shown in \cref{sec:HC}.
Let $X_{1},\dots,X_{n}$ be i.i.d. random variables taking values in a measurable space $(\mathcal{X},\mathcal{A})$ with a common distribution $P$ and $n\geq3$.
Let $\mathcal{F}$ be a pointwise measurable class of measurable functions $\mathcal{X}\to\mathbb{R}$ with an envelope function $F$ such that $Pf=0$ for every $f\in\mathcal{F}$.
Let $Z=\sup\{\mathbb{G}_{n}f:f\in\mathcal{F}\}$, where $\mathbb{G}_{n}=n^{-1/2}\sum_{i=1}^{n}(\delta_{X_{i}}-P)$ is the empirical process associated to $X_{1},\dots,X_{n}$.
Let $G_{P}$ be a tight Gaussian process in $\ell^{\infty}(\mathcal{F})$ with mean zero and covariance function $\mathbb{E}(G_{P}fG_{P}g)=P(fg)$ for every $f,g\in\mathcal{F}$.

\begin{theorem}
\label{thm:CCK14}
Suppose that there are constants $v,b\in(0,\infty)$, $q\in[4,\infty]$, and $\delta_{3},\delta_{4}\in(0,1)$ such that $\vc(\mathcal{F})\leq v$, $\|F\|_{P,q}\leq b$, and $\sup\{\|f\|_{P,k}:f\in\mathcal{F}\}\leq\delta_{k}\|F\|_{P,k}$ for $k=3,4$.
Then for every $\gamma\in(0,1)$, there is a random variable $\tilde{Z}=_{d}\sup\{G_{P}f:f\in\mathcal{F}\}$ such that
\begin{equation}
\mathbb{P}(|Z-\tilde{Z}|>\Delta)\leq C_{q}\Big(\gamma+\frac{\log n}{n}\Big),
\end{equation}
where
\begin{equation}
\Delta=\frac{bK}{\gamma^{1/2}n^{1/2-1/q}}+\frac{b\delta_{4}K^{3/4}}{\gamma^{1/2}n^{1/4}}+\frac{b\delta_{3}K^{2/3}}{\gamma^{1/3}n^{1/6}},
\label{eq:Delta}
\end{equation}
$K=c_{q}v\log(An/\delta_{3}\delta_{4})$, $c_{q}$ and $C_{q}$ are positive constants depending only on $q$, and $A>0$ is an absolute constant.
\end{theorem}

\begin{remark}
The existence of a tight choice of $G_{P}$ having, with probability one, uniformly $L_{2}(P)$-continuous sample functions follows from the same argument as that in proving Lemma 2.1 in \cite{CCK14-MR3262461}, which is based on $\vc(\mathcal{F})\leq v$, $\|F\|_{P,q}\leq b$, and an application of Example 1.5.10 in \cite{VW23-MR4628026} and Theorem 1.1 in \cite{Dud73-MR346884}.
\end{remark}
%
%

\iffalse
\begin{remark}
We note that \cref{thm:CCK14} is a slight generalization of Corollary 2.2 in \cite{CCK14-MR3262461}.
By a careful examination of their proof, when $\delta_{3}=(\sigma/b)^{2/3}$ and $\delta_{4}=(\sigma/b)^{1/2}$ with $\sigma=\sup\{\|f\|_{P,2}:f\in\mathcal{F}\}\in(0,b]$, it follows that $K=c_{q}v\log(Abn/\sigma)$ and that
\begin{equation}
\label{eq:Delta-CCK14}
\Delta=\frac{bK}{\gamma^{1/2}n^{1/2-1/q}}+\frac{(b\sigma)^{1/2}K^{3/4}}{\gamma^{1/2}n^{1/4}}+\frac{(b\sigma^{2})^{1/3}K^{2/3}}{\gamma^{1/3}n^{1/6}},
\end{equation}
which recovers their original result.
However, it can be shown that
\begin{equation}
\label{eq:lem-CCK-counterexample}
\frac{\sup_{f\in\mathcal{F}}\|f\|_{P,3}}{(\sigma/b)^{2/3}\|F\|_{P,3}}\to\infty\quad\text{and}\quad\frac{\sup_{f\in\mathcal{F}}\|f\|_{P,4}}{(\sigma/b)^{1/2}\|F\|_{P,4}}\to\infty
\end{equation}
for the multi-level thresholding statistic \eqref{eq:multi-level-thresholding-statistic-original}.
By definition, we cannot take $\delta_{3}=(\sigma/b)^{2/3}$ and $\delta_{4}=(\sigma/b)^{1/2}$ in this case, which means that the error bound \eqref{eq:Delta-CCK14} cannot be directly applied.
This is the motivation to provide a generalization.
\end{remark}
\fi
%
%

For completeness, we provide the proof of \cref{thm:CCK14} in \cref{pf:lem-CCK14} although the argument is quite similar to the mentioned paper.
Next, we extend \cref{thm:CCK14} for dependent observations by invoking the big-block-small-block technique.
Let $\beta(\mathcal{A},\mathcal{B})=\mathbb{E}\ess\sup\{|\mathbb{P}(B|\mathcal{A})-\mathbb{P}(B)|:B\in\mathcal{B}\}$ be the $\beta$-mixing coefficient between $\sigma$-algebras $\mathcal{A}$ and $\mathcal{B}$.
Let $X_{1},\dots,X_{n}$ be stationary random variables taking values in a Polish space $\mathcal{X}$ with a common distribution $P$ and $\beta$-mixing coefficients $\beta(\ell)=\sup\{\beta(\sigma(X_{i}:1\leq i\leq k),\sigma(X_{i}:i\geq k+\ell)):k\geq1\}$ for $\ell\geq1$.
Let $a_{1}\geq a_{2}\geq1$ be block sizes and $m=[n/(a_{1}+a_{2})]\geq3$ be the number of blocks.
Let $\mathcal{F}$ be a pointwise measurable class of measurable functions $\mathcal{X}\to\mathbb{R}$ with an envelope function $F$ such that $Pf=0$ for every $f\in\mathcal{F}$.
Let $Z=\sup\{\mathbb{G}_{n}f:f\in\mathcal{F}\}$, where $\mathbb{G}_{n}=n^{-1/2}\sum_{i=1}^{n}(\delta_{X_{i}}-P)$ is the empirical process associated to $X_{1},\dots,X_{n}$.
Let $G_{P_{0}}$ be a tight Gaussian process in $\ell^{\infty}(\mathcal{F})$ with mean zero and covariance function $\sigma_{P_{0}}(f,g)=\mathbb{E}(\mathbb{G}_{n}f\mathbb{G}_{n}g)$ for every $f,g\in\mathcal{F}$.
Let $\sigma_{P_{1}}(f,g)=(ma_{1}/n)\mathbb{E}(\mathbb{G}_{a_{1}}f\mathbb{G}_{a_{1}}g)$.

\begin{theorem}%[Gaussian approximation for empirical processes indexed by VC classes of functions with stationary $\beta$-mixing observations]
\label{thm:GA-mixing}
Suppose that there are constants $v,b\in(0,\infty)$, $q\in[4,\infty]$, and $\delta_{2},\delta_{3},\delta_{4}\in(0,1)$ such that $\vc(\mathcal{F})\leq v$, $\|F\|_{P,q}\leq b$, and $\sup\{\|f\|_{P,k}:f\in\mathcal{F}\}\leq\delta_{k}\cdot P\text{-}\ess\inf F$ for $k=2,3,4$.
Then for every $\gamma\in(0,1)$, there is a random variable $\tilde{Z}_{0}=_{d}\sup\{G_{P_{0}}f:f\in\mathcal{F}\}$ such that
\begin{equation}
\mathbb{P}(|Z-\tilde{Z}_{0}|>\Delta_{0})\leq C_{q}\Big(\gamma+\frac{\log m}{m}\Big)+2m\beta(a_{2}),
\end{equation}
where
\begin{align}
\Delta_{0}&=\frac{bKn^{1/2}}{\gamma m^{1-1/q}}+\frac{b\delta_{4}K^{3/4}n^{1/2}}{\gamma^{1/2}m^{3/4}}+\frac{b\delta_{3}K^{2/3}n^{1/2}}{\gamma^{1/3}m^{2/3}}\notag\\&\quad+\frac{K^{1/2}}{\gamma^{1/2}}\|\sigma_{P_{0}}-\sigma_{P_{1}}\|_{\mathcal{F}^{2}}^{1/2}+\frac{b\delta_{2}K^{1/2}n^{1/2}}{(a_{1}/a_{2})\gamma m^{1/2}},
\label{eq:Delta-0}
\end{align}
$K=c_{q}v\log(Am/\delta_{2}\delta_{3}\delta_{4})$, $c_{q}$ and $C_{q}$ are positive constants depending only on $q$, and $A>0$ is an absolute constant.
\end{theorem}

\begin{remark}
\label{rmk:thm-GA-mixing}
When $\beta(a_{2})\leq\exp(-c_{1}a_{2}^{c_{2}})$ for some constants $c_{1},c_{2}>0$, we may take $a_{1}=(\log n)^{\tau_{1}}$ and $a_{2}=(\log n)^{\tau_{2}}$ for $\tau_{1}\geq \tau_{2}>1/c_{2}$, and then conclude that $m\beta(a_{2})\to0$ as $n\to\infty$ and that \eqref{eq:Delta-0} becomes
\begin{align}
\Delta_{0}&=\frac{bK(\log n)^{(1-1/q)\tau_{1}}}{\gamma n^{1/2-1/q}}+\frac{b\delta_{4}K^{3/4}(\log n)^{3\tau_{1}/4}}{\gamma^{1/2}n^{1/4}}+\frac{b\delta_{3}K^{2/3}(\log n)^{2\tau_{1}/3}}{\gamma^{1/3}n^{1/6}}\notag\\&\quad+\frac{K^{1/2}}{\gamma^{1/2}}\|\sigma_{P_{0}}-\sigma_{P_{1}}\|_{\mathcal{F}^{2}}^{1/2}+\frac{b\delta_{2}K^{1/2}}{\gamma(\log n)^{\tau_{1}/2-\tau_{2}}}.\label{eq:Delta-0-log}
\end{align}
The first three terms in \eqref{eq:Delta-0-log} are almost the same as those in \eqref{eq:Delta} up to some polylogarithmic factors while the last two terms stem from the application of the big-block-small-block technique due to data dependency.
\end{remark}

\subsection{Proof sketch of \cref{thm:GA-mixing}}
We defer the rigorous proof of \cref{thm:GA-mixing} to \cref{pf:thm-GA-mixing} and sketch it below.
The key to prove \cref{thm:GA-mixing} via the big-block-small-block technique is to apply Berbee's lemma (see \cref{lem:Berbee} in \cref{pf:thm-GA-mixing}) for the stationary $\beta$-mixing observations.
Then the original empirical process $\mathbb{G}_{n}$ can be approximated by the empirical process associated to i.i.d. big blocks so that \cref{thm:CCK14} can be applied with a slight modification (see \cref{lem:GA-blocking}).
The last terms involving $b\delta_{2}$ in \eqref{eq:Delta-0} and \eqref{eq:Delta-0-log} come from the approximation error by dropping the small blocks.

More specifically, let $P_{0}$, $P_{1}$, and $P_{2}$ be the distributions of $(X_{1},\dots,X_{n})$, $(X_{1},\dots,X_{a_{1}})$, and $(X_{1},\dots,X_{a_{2}})$, respectively.
Define the blocks as
\begin{align}
Y_{2i-1}&=(X_{(a_{1}+a_{2})(i-1)+1},\dots,X_{(a_{1}+a_{2})(i-1)+a_{1}}),&1\leq i\leq m,\\Y_{2i}&=(X_{(a_{1}+a_{2})(i-1)+a_{1}+1},\dots,X_{(a_{1}+a_{2})i}),&1\leq i\leq m,\\Y_{2m+1}&=(X_{(a_{1}+a_{2})m+1},\dots,X_{n}).&
\end{align}
Let $\varphi_{j}:\mathbb{R}^{\mathcal{X}}\to\mathbb{R}^{\mathcal{X}^{a_{j}}}$ with $f\mapsto(m/n)^{1/2}\sum_{i=1}^{a_{j}}f\circ\pi_{i}$ for $j\in\{1,2\}$, where $\pi_{i}(x)=x_{i}$ denotes the projection onto the $i$-th coordinate for a vector $x=(x_{1},\dots,x_{n})$.
For instance, $\varphi_{1}(f)(Y_{2i-1})=(m/n)^{1/2}\sum_{j=1}^{a_{1}}f(X_{(a_{1}+a_{2})(i-1)+j})$.
Observe that we may decompose $\mathbb{G}_{n}$ as
\begin{equation}
\label{eq:Gn-decomposition}
\mathbb{G}_{n}f=\mathbb{G}_{n,1}(\varphi_{1}f)+\mathbb{G}_{n,2}(\varphi_{2}f)+\mathbb{G}_{n,3}f,
\end{equation}
where $\mathbb{G}_{n,1}=m^{-1/2}\sum_{i=1}^{m}(\delta_{Y_{2i-1}}-P_{1})$ is the empirical process associated to big blocks $Y_{1},\dots,Y_{2m-1}$, $\mathbb{G}_{n,2}=m^{-1/2}\sum_{i=1}^{m}(\delta_{Y_{2i}}-P_{2})$ is the empirical process associated to small blocks $Y_{2},\dots,Y_{2m}$, and $\mathbb{G}_{n,3}$ is the remaining fraction of $\mathbb{G}_{n}$ associated to $Y_{2m+1}$.
Let $\mathbb{G}_{n,1}^{*}=m^{-1/2}\sum_{i=1}^{m}(\delta_{Y_{2i-1}^{*}}-P_{1})$ be the empirical process associated to $Y_{1}^{*},\dots,Y_{2m-1}^{*}$, the i.i.d. coupling version of $Y_{1},\dots,Y_{2m-1}$.
The general idea is to approximate $\mathbb{G}_{n}$ by $\mathbb{G}_{n,1}(\varphi_{1}\cdot)$, and then by $\mathbb{G}_{n,1}^{*}(\varphi_{1}\cdot)$.

Observe that $\sigma_{P_{0}}$ is the covariance function of $\mathbb{G}_{n}$ while $\sigma_{P_{1}}$ is the covariance function of $\mathbb{G}_{n,1}^{*}(\varphi_{1}\cdot)$.
To control the approximation error between the natural limiting Gaussian processes of the two empirical processes, we shall invoke the following coupling inequality for comparisons between suprema of Gaussian processes whose proof can be found in \cref{pf:lem-comparison-process}.
Note that the lemma eventually produces the terms bearing $\|\sigma_{P_{0}}-\sigma_{P_{1}}\|_{\mathcal{F}^{2}}$ in \eqref{eq:Delta-0} and \eqref{eq:Delta-0-log}.
Let $N(\varepsilon,\mathcal{F},L_{r}(Q))$ be the $\varepsilon$-covering number of $\mathcal{F}$ related to the $L_{r}(Q)$-seminorm, namely the minimum number of balls $\{g:\|g-f\|_{Q,r}<\varepsilon\}$ of radius $\varepsilon$ needed to cover $\mathcal{F}$, where $\varepsilon>0$ and $r\in[1,\infty]$.
Let $\mathcal{F}_{Q,\varepsilon}=\{f-g:f,g\in\mathcal{F},\|f-g\|_{Q,2}<\varepsilon\|F\|_{Q,2}\}$.

\begin{lemma}%[Coupling inequality for comparisons between Gaussian processes indexed by classes of functions under different transformations]
\label{lem:comparison-process}
Let $P$ and $Q$ be two distributions on measurable spaces $(\mathcal{Y},\mathcal{A})$ and $(\mathcal{Z},\mathcal{B})$, respectively.
Let $\mathcal{F}$ be a pointwise measurable class of measurable functions $\mathcal{X}\to\mathbb{R}$.
Let $\varphi_{P}:\mathbb{R}^{\mathcal{X}}\to\mathbb{R}^{\mathcal{Y}}$ and $\varphi_{Q}:\mathbb{R}^{\mathcal{X}}\to\mathbb{R}^{\mathcal{Z}}$ be two transformations on $\mathcal{F}$.
For each $R\in\{P,Q\}$, let $F_{R}$ be the envelope function of $\varphi_{R}\mathcal{F}$.
Suppose that $R(\varphi_{R}f)=0$ for every $f\in\mathcal{F}$ and $\|F_{R}\|_{R,2}<\infty$.
Let $G_{R}$ be a tight Gaussian process in $\ell^{\infty}(\varphi_{R}\mathcal{F})$ with mean zero and covariance function $\sigma_{R}(f,g)=\mathbb{E}(G_{R}(\varphi_{R}f)G_{R}(\varphi_{R}g))=R(\varphi_{R}f\varphi_{R}g)$ for every $f,g\in\mathcal{F}$.
Let $Z_{P}=\sup\{G_{P}(\varphi_{P}f):f\in\mathcal{F}\}$.
Then for every $\beta>0$, $\gamma\in(0,1)$, $\delta>1/\beta$, and $\varepsilon\in(0,1)$, there is a random variable $\tilde{Z}_{Q}=_{d}\sup\{G_{Q}(\varphi_{Q}f):f\in\mathcal{F}\}$ such that
\begin{equation}
\mathbb{P}(|Z_{P}-\tilde{Z}_{Q}|>\Delta_{P,Q}(\beta,\gamma,\delta,\varepsilon))\leq\gamma+\frac{\eta+C\beta\delta^{-1}\|\sigma_{P}-\sigma_{Q}\|_{\mathcal{F}^{2}}}{1-\eta},
\end{equation}
where
\begin{align}
\Delta_{P,Q}(\beta,\gamma,\delta,\varepsilon)&=2\beta^{-1}\log N_{P,Q}(\varepsilon)+3\delta+b_{P,Q}(\gamma,\varepsilon),\\N_{P,Q}(\varepsilon)&={\textstyle \sum_{R\in\{P,Q\}}}N(\varepsilon\|F_{R}\|_{R,2},\varphi_{R}\mathcal{F},L_{2}(R)),\\b_{P,Q}(\gamma,\varepsilon)&={\textstyle \sum_{R\in\{P,Q\}}}(\mathbb{E}\|G_{R}\|_{(\varphi_{R}\mathcal{F})_{R,\varepsilon}}+\varepsilon\|F_{R}\|_{R,2}(2\log2/\gamma)^{1/2}),
\end{align}
$\eta=e^{-\alpha/2}(1+\alpha)^{1/2}<1$, $\alpha=\beta^{2}\delta^{2}-1>0$, and $C$ is an absolute constant.
\end{lemma}

\section{Higher criticism asymptotics}
\label{sec:HC}
We now consider the higher criticism statistic $\HC_{p}^{*}$ in \eqref{eq:HC-significance-level} when $T_{j}$'s are stationary $N(0,1)$ random variables with $\beta$-mixing coefficients $\beta(\cdot)$, and autocorrelation coefficients $\rho_{k}=\mathbb{E}(T_{j}T_{j+k})$.
Let $\pi_{0}(\lambda)=\mathbb{P}(|N(0,1)|\geq\lambda)$ and $\sigma_{0}^{2}(\lambda)=\pi_{0}(\lambda)(1-\pi_{0}(\lambda))$ for $\lambda\geq0$.
Then the $j$-th individual $p$-value is given by $p_{j}=\pi_{0}(|T_{j}|)$.
Let $[\alpha_{1},\alpha_{2}]=[(\log p)^{c}/p,1/(\log p)^{d}]$ be the searching range of candidate significance levels, where $c,d>0$ are two constants.
Let $\lambda_{1}$ and $\lambda_{2}$ be two quantities satisfying $\pi_{0}(\lambda_{1})=\alpha_{2}$ and $\pi_{0}(\lambda_{2})=\alpha_{1}$.
Let $\mathbb{G}_{p}=p^{-1/2}\sum_{j=1}^{p}(\delta_{T_{j}}-\Phi)$ be the empirical process associated to $T_{1},\dots,T_{p}$, where $\Phi$ is the distribution of $N(0,1)$.
To apply the Gaussian approximation result in \cref{thm:GA-mixing}, we rewrite \eqref{eq:HC-significance-level} as
\begin{equation}
\label{eq:HC-empirical-process}
\HC_{p}^{*}=\sup_{\lambda\in[\lambda_{1},\lambda_{2}]}\mathbb{G}_{p}f_{\lambda},\quad\text{where}\quad f_{\lambda}(x)=\frac{1}{\sigma_{0}(\lambda)}(1_{[\lambda,\infty)}(|x|)-\pi_{0}(\lambda)).
\end{equation}
Clearly, $\Phi f_{\lambda}=0$ and $\Phi f_{\lambda}^{2}=1$ for every $\lambda\in[\lambda_{1},\lambda_{2}]$.
Let $\mathcal{F}=\{f_{\lambda}:\lambda\in[\lambda_{1},\lambda_{2}]\}$.
Then $F(x)=1/\sigma_{0}(\lambda_{2})$ is an envelope function of $\mathcal{F}$.
Let $\mathbb{B}_{p}$ be a tight Gaussian process in $\ell^{\infty}(\mathcal{F})$ with mean zero and covariance function $\mathbb{E}(\mathbb{B}_{p}f_{\lambda}\mathbb{B}_{p}f_{\nu})=\mathbb{E}(\mathbb{G}_{p}f_{\lambda}\mathbb{G}_{p}f_{\nu})$.
The following theorem provides a Gaussian approximation result for higher criticism with dependent test statistics.

\begin{theorem}%[Gaussian approximation for higher criticism with dependent observations]
\label{thm:GA-HC}
Suppose that there are constants $c_{1},c_{2}\in(0,\infty)$ and $\rho_{0}\in[0,1)$ such that $\beta(\ell)\leq\exp(-c_{1}\ell^{c_{2}})$ for $\ell\geq1$, $\sum_{k=1}^{\infty}|\rho_{k}|<\infty$, and $\sup\{|\rho_{k}|:k\geq1\}\leq\rho_{0}$.
Let $c>8+8/c_{2}$ and $d>3(1+\rho_{0})/2(1-\rho_{0})$.
Then there is a constant $\tau>0$ and a random variable $\HC_{p}^{**}=_{d}\sup\{\mathbb{B}_{p}f_{\lambda}:\lambda\in[\lambda_{1},\lambda_{2}]\}$ such that
\begin{equation}
\label{eq:thm-GA-HC}
|\HC_{p}^{*}-\HC_{p}^{**}|=O_{\mathbb{P}}((\log p)^{-\tau}).
\end{equation}
\end{theorem}

\begin{remark}
We do not try to optimize the approximation error bound in \eqref{eq:thm-GA-HC} as the polylogarithmic rate is fast enough to derive the asymptotic distribution of the higher criticism statistic as will be shown in the next theorem.
\end{remark}

Let $\mathbb{B}_{p}^{0}$ be a Gaussian process modified from $\mathbb{B}_{p}$ as if the test statistics $T_{j}$ were independent $N(0,1)$ variables such that $\mathbb{B}_{p}^{0}f_{\lambda}=B^{0}(\pi_{0}(\lambda))/\sigma_{0}(\lambda)$ for every $\lambda\in[\lambda_{1},\lambda_{2}]$, where $B^{0}$ is the standard Brownian bridge process in \cref{sec:intro}.
As shown in the following theorem, an attracting feature of the higher criticism is that the dependence among test statistics can be negligible by searching large candidate thresholding levels $\lambda$, or equivalently searching small candidate significance levels $\alpha$.
The covariance function of the limiting Gaussian process $\mathbb{B}_{p}$ is almost the same as that for the independent case $\mathbb{B}_{p}^{0}$, which leads to the asymptotic distribution of the higher criticism statistic under dependence.

\begin{theorem}
\label{thm:HC-limit}
Let the assumptions in \cref{thm:GA-HC} be valid.
Then
\begin{equation}
\label{eq:HC-almost-independent}
\sup_{\lambda,\nu\in[\lambda_{1},\lambda_{2}]}\Big|\mathbb{E}(\mathbb{B}_{p}f_{\lambda}\mathbb{B}_{p}f_{\nu})-\mathbb{E}(\mathbb{B}_{p}^{0}f_{\lambda}\mathbb{B}_{p}^{0}f_{\nu})\Big|=O\Big(\frac{(\log\log p)^{(1-\rho_{0})/2(1+\rho_{0})}}{(\log p)^{(1-\rho_{0})d/(1+\rho_{0})-1/2}}\Big).
\end{equation}
Consequently, the approximation error bound in \eqref{eq:HC-significance-level-GA} and the asymptotic distribution in \eqref{eq:HC-significance-level-asymptotics} are valid for the higher criticism in this case.
\end{theorem}

As the same asymptotic distribution of the dependent higher criticism is the same as that for the independent case, it is safe to conclude that the higher criticism is robust to the dependence among different individual testing problems.
As shown in \cref{sec:DB}, the range of $[\alpha_{1},\alpha_{2}]=[(\log p)^{c}/p,1/(\log p)^{d}]$ is enough for the higher criticism test to attain the optimal detection boundary.

\subsection{Higher criticism $t$-statistic}
\label{sec:HC-t}
We next consider the case in which the test statistics $T_{j}$ are not necessarily $N(0,1)$ variables (but still stationary $\beta$-mixing).
Instead, we assume that $T_{j}=n^{1/2}\bar{Y}_{j}/\hat{\sigma}_{j}$ is the $j$-th individual $t$-statistic, where $\bar{Y}_{j}=n^{-1}\sum_{i=1}^{n}Y_{ij}$, $\hat{\sigma}_{j}^{2}=(n-1)^{-1}\sum_{i=1}^{n}(Y_{ij}-\bar{Y}_{j})^{2}$, and $Y_{1},\dots,Y_{n}$ are i.i.d. $\mathbb{R}^{p}$-valued random vectors with mean zero and $\sup\{\mathbb{E}|Y_{ij}|^{2+\delta}:j\geq1\}<\infty$ for $\delta\in(0,1]$.
Let $p_{j}=\pi_{0}(|T_{j}|)$ be the $j$-th individual asymptotic $p$-value.
Let the higher criticism $t$-statistic be defined as in \eqref{eq:HC-significance-level} with $[\alpha_{1},\alpha_{2}]=[(\log p)^{c}/p,1/p^{d}]$ for two constants $c,d>0$.
Let $\lambda_{1}$ and $\lambda_{2}$ be two quantities satisfying $\pi_{0}(\lambda_{1})=\alpha_{2}$ and $\pi_{0}(\lambda_{2})=\alpha_{1}$.
The following theorem provides a Gaussian approximation result for the higher criticism $t$-statistic.

\begin{theorem}
\label{thm:GA-HC-t}
Suppose that there are constants $c_{1},c_{2}\in(0,\infty)$ and $\rho_{0}\in[0,1)$ such that $\beta(\ell)\leq\exp(-c_{1}\ell^{c_{2}})$ for $\ell\geq1$, $\sum_{k=1}^{\infty}|\rho_{k}|<\infty$, and $\sup\{|\rho_{k}|:k\geq1\}\leq\rho_{0}$.
Let $c>8+8/c_{2}$, $d\in(\rho_{0}^{2},1)$, and $\log p=o(n^{\delta/(4+2\delta)})$.
Then there is a constant $\tau>0$ and a random variable $\BB_{p}^{*}$ described as in \cref{sec:intro} such that
\begin{equation}
\label{eq:thm-GA-HC-t-1}
|\HC_{p}^{*}-\BB_{p}^{*}|=O_{\mathbb{P}}\Big(\frac{1}{(\log p)^{\tau}}+\frac{p^{(1-d)/2}(\log p)^{1+\delta/2}}{n^{\delta/2}}\Big).
\end{equation}
In particular, when $p=n^{1/\theta}$ and $d>(1-\delta\theta)\vee\rho_{0}^{2}$,
\begin{equation}
\mathbb{P}\Big(\HC_{p}^{*}\leq\frac{2\log\log p+(1/2)\log\log\log p+x}{(2\log\log p)^{1/2}}\Big)\to\exp\Big(-\frac{1-d}{2\pi^{1/2}}e^{-x}\Big).
\end{equation}
\end{theorem}

When $\delta=1$ and $\rho_{0}=0$, the searching range $[\alpha_{1},\alpha_{2}]=[(\log p)^{c}/p,1/p^{d}]$ can approach to that in \cite{DHJ11-MR2815777} for independent $t$-statistics up to some polylogarithmic terms.
The condition of $d>1-\delta\theta$ is necessary to control the second term in \eqref{eq:thm-GA-HC-t-1} introduced due to the non-Gaussianity, which disappears in \eqref{eq:thm-GA-HC} for the Gaussian case.
Note that \cref{thm:GA-HC-t} additionally provides the asymptotic distribution of the higher criticism $t$-statistic under dependence that was not covered by \cite{HJ10-MR2662357} or \cite{DHJ11-MR2815777}.
On the other hand, we do not require the sub-Gaussian condition used in \cite{ZCX13-MR3161449} and \cite{CLZ19-MR3911118} and do not need to omit the range of $\alpha\in[(\log p)^{c}/p,L_{p}/p^{1-\eta}]$ that is crucial for highly sparse signal detection as discussed in \cref{sec:DB}.
By invoking the following bivariate Cram\'er type moderate deviation result as a direct extension of Lemma 8.1 in \cite{QCS25+}, we reduce the exponential type moment condition to a finite ($2+\delta$)-th moment condition, which is more suitable for heavy-tailed data.
Note that the condition of $d>\rho_{0}^{2}$ is required for the bivariate result \eqref{eq:MDR-2}.

\begin{lemma}
\label{lem:MDR}
Let $(Z_{j},Z_{k})\sim N(0,(\begin{smallmatrix}1 & \rho_{|j-k|}\\
\rho_{|j-k|} & 1
\end{smallmatrix}))$.
Then we have
\begin{equation}
\label{eq:MDR-1}
\sup_{j\geq1}\sup_{\lambda\in[\lambda_{1},\lambda_{2}]}\Big|\frac{\mathbb{P}(|T_{j}|\geq\lambda)}{\mathbb{P}(|Z_{j}|\geq\lambda)}-1\Big|=O\Big(\frac{(\log p)^{1+\delta/2}}{n^{\delta/2}}\Big),
\end{equation}
and
\begin{equation}
\label{eq:MDR-2}
\sup_{j,k\geq1}\sup_{\lambda,\nu\in[\lambda_{1},\lambda_{2}]}\Big|\frac{\mathbb{P}(|T_{j}|\geq\lambda,|T_{k}|\geq\nu)}{\mathbb{P}(|Z_{j}|\geq\lambda,|Z_{k}|\geq\nu)}-1\Big|=O\Big(\frac{(\log p)^{1+\delta/2}}{n^{\delta/2}}\Big).
\end{equation}
\end{lemma}

\subsection{Extension for multi-level thresholding $t$-statistic}
\label{sec:MTT}
Although the functions $f_{\lambda}$ indexed in the empirical process associated to the higher criticism are the normalized indicator functions, the Gaussian approximation results in Theorems \ref{thm:CCK14} and \ref{thm:GA-mixing} allow for more complex classes of functions.
The asymptotic results for the higher criticism $t$-statistic can be further extended for the more general multi-level thresholding $t$-statistic.
When $T_{j}$'s are dependent $z$-statistics, the multi-level thresholding test has been discussed in \cite{ZCX13-MR3161449}.
In particular, Theorem 4 in the mentioned paper shows that the test is asymptotically more powerful than the higher criticism test under some fixed alternatives.

Let $Y_{1},\dots,Y_{n}$ be i.i.d. $\mathbb{R}^{p}$-valued random vectors defined as before.
Let $[\lambda_{1},\lambda_{2}]\subset(0,\infty)$ be the searching range of candidate thresholding levels and $M>0$ be a truncating level.
The multi-level thresholding $t$-statistic is given by
\begin{equation}
\label{eq:MT-threhsolding-t}
\MT_{p}^{*}=\sup_{\lambda\in[\lambda_{1},\lambda_{2}]}\frac{1}{p^{1/2}\sigma_{0}(\lambda)}\sum_{j=1}^{p}(|T_{j}|^{2}1_{[\lambda,M]}(|T_{j}|)-\mu_{0}(\lambda)),
\end{equation}
where
\begin{equation}
\label{eq:multi-level-thresholding-statistic-mean}
\mu_{0}(\lambda)=2(\lambda\phi(\lambda)-M\phi(M)+\Phi(M)-\Phi(\lambda)),
\end{equation}
and
\begin{equation}
\label{eq:multi-level-thresholding-statistic-variance}
\sigma_{0}^{2}(\lambda)=2(\lambda^{3}\phi(\lambda)-M^{3}\phi(M))+(3-\mu_{0}(\lambda))\mu_{0}(\lambda).
\end{equation}
Here, $\lambda_{1}$ and $\lambda_{2}$ are two quantities satisfying $\sigma_{0}^{2}(\lambda_{1})=\alpha_{2}$ and $\sigma_{0}^{2}(\lambda_{2})=\alpha_{1}$, where $[\alpha_{1},\alpha_{2}]=[(\log p)^{c}/p,1/p^{d}]$ for two constants $c,d>0$.

Write $\MT_{p}^{*}=\sup\{\MT_{p,\lambda}:\lambda\in[\lambda_{1},\lambda_{2}]\}$.
We notice that \cite{QCS25+} has proven the asymptotic normality of the single-level thresholding $t$-statistic $\MT_{p,\lambda}$ for each candidate thresholding level $\lambda\in[\lambda_{1},\lambda_{2}]$ by using the bivariate Cram\'er type moderate deviation result in \cref{lem:MDR}.
The following theorem further provides a Gaussian approximation result for the multi-level thresholding $t$-statistic.
We emphasize that the natural limiting Gaussian process in this case is no longer based on the standard Brownian bridge $B^{0}$ for the higher criticism setting.
Instead, we shall define a tight Gaussian process $\mathbb{B}_{p}$ indexed by $[\lambda_{1},\lambda_{2}]$ with mean zero and covariance function
\begin{equation}
\label{eq:cov-brownian-bridge-MT}
\mathbb{E}(\mathbb{B}_{p,\lambda}\mathbb{B}_{p,\nu})=\frac{\sigma_{0}(\nu)}{\sigma_{0}(\lambda)}+\frac{\mu_{0}(\nu)}{\sigma_{0}(\lambda)\sigma_{0}(\nu)}(\mu_{0}(\nu)-\mu_{0}(\lambda)),
\end{equation}
where $\lambda_{1}\leq\lambda\leq\nu\leq\lambda_{2}$.

\begin{theorem}
\label{thm:GA-MT-t}
Suppose that there are constants $c_{1},c_{2}\in(0,\infty)$ and $\rho_{0}\in[0,1)$ such that $\beta(\ell)\leq\exp(-c_{1}\ell^{c_{2}})$ for $\ell\geq1$, $\sum_{k=1}^{\infty}|\rho_{k}|<\infty$, and $\sup\{|\rho_{k}|:k\geq1\}\leq\rho_{0}$.
Let $c>8+8/c_{2}$, $d\in(\rho_{0}^{2},1)$, $\log p=o(n^{\delta/(4+2\delta)})$, and $M=(2\log p)^{1/2}$.
Then there is a constant $\tau>0$ and a random variable $\mathbb{B}_{p}^{*}=_{d}\sup\{\mathbb{B}_{p,\lambda}:\lambda\in[\lambda_{1},\lambda_{2}]\}$ such that
\begin{equation}
\label{eq:thm-GA-MT}
|\MT_{p}^{*}-\mathbb{B}_{p}^{*}|=O_{\mathbb{P}}\Big(\frac{1}{(\log p)^{\tau}}+\frac{p^{(1-d)/2}(\log p)^{\delta/2}}{n^{\delta/2}}\Big).
\end{equation}
In particular, when $p=n^{1/\theta}$ and $d>(1-\delta\theta)\vee\rho_{0}^{2}$,
\begin{equation}
\label{eq:MT-limit}
\mathbb{P}\Big(\MT_{p}^{*}\leq\frac{2\log\log p+(1/2)\log\log\log p+x}{(2\log\log p)^{1/2}}\Big)\to\exp\Big(-\frac{1-d}{2\pi^{1/2}}e^{-x}\Big).
\end{equation}
\end{theorem}

To address the asymptotic distribution of $\mathbb{B}_{p}^{*}$ that differs from $\BB_{p}^{*}$, we shall invoke the following limit result, which is an immediate generalization of Theorem 4.2 in \citet{Hus90-MR1062062} for sequences of asymptotically locally stationary Gaussian processes.
The proof of the lemma follows the same lines in the mentioned paper and we omit it.

\begin{lemma}
\label{lem:Hus90-modification}
Let $\{\mathbb{X}_{n}(t):t\geq0\}$ be a sequence of real Gaussian processes with mean zero, variance one, continuous sample functions, and the covariance function satisfying that as $h\to0$ and $n\to\infty$,
\begin{equation}
\label{eq:asymptotically-locally-stationary-Gaussian-process}
\mathbb{E}(\mathbb{X}_{n}(t+h)\mathbb{X}_{n}(t))=1-C(t)|h|^{\alpha}+o(|h|^{\alpha})
\end{equation}
uniformly in $t\geq0$ for a constant $\alpha\in(0,2]$ and a uniformly continuous function $C$ with $0<\inf\{C(t):t\geq0\}\leq\sup\{C(t):t\geq0\}<\infty$.

Suppose in addition that $\delta_{n}(\tau)\log\tau\to0$ as $\tau\to\infty$ and $n\to\infty$, and that $\mathbb{E}(\mathbb{X}_{n}(s)\mathbb{X}_{n}(t))=1$ if and only if $s=t$, where $\delta_{n}(\tau)=\sup\{\mathbb{E}(\mathbb{X}_{n}(s)\mathbb{X}_{n}(t)):|s-t|\geq\tau\}$.
Then as $T\to\infty$ and $n\to\infty$,
\begin{equation}
\mathbb{P}\Big( a_{T}\Big(\sup_{t\in[0,T]}\mathbb{X}_{n}(t)-b_{T}\Big)\leq x\Big)\to\exp(-e^{-x}),
\end{equation}
where $a_{T}=(2\log T)^{1/2}$,
\begin{equation}
b_{T}=a_{T}+\frac{1}{a_{T}}\Big(\frac{2-\alpha}{2\alpha}\log\log T+\log(C_{T}^{*}H_{\alpha}(2\pi)^{-1/2}2^{(2-\alpha)/2\alpha})\Big),
\end{equation}
$C_{T}^{*}=T^{-1}\int_{0}^{T}C^{1/\alpha}(t)dt$, and $H_{\alpha}$ is a positive constant depending only on the index $\alpha$ with $H_{1}=1$ and $H_{2}=\pi^{-1/2}$; cf. Lemmas 12.2.7 and 12.2.8 in \cite{LLR83-MR691492}.
\end{lemma}

\section{Proof of Theorems \ref{thm:GA-HC} and \ref{thm:HC-limit}}
\label{pf:thm-GA-HC}
The key of the proof is to determine the quantities in \cref{thm:GA-mixing} for the higher criticism class $\mathcal{F}=\{f_{\lambda}:\lambda\in[\lambda_{1},\lambda_{2}]\}$.
By the construction in \eqref{eq:HC-empirical-process}, $b=1/\sigma_{0}(\lambda_{2})$ and $q=\infty$.
Then it suffices to determine $v$, $\delta_{2}$, $\delta_{3}$, $\delta_{4}$, and $\|\sigma_{P_{0}}-\sigma_{P_{1}}\|_{\mathcal{F}^{2}}$.

\subsubsection*{Step 1: Determining $v$}
We shall first verify that $\mathcal{F}$ is a VC-subgraph class.
In view of Example 2.6.24 in \cite{VW23-MR4628026}, the class $\{1_{[\lambda,\infty)}(|x|)/\sigma_{0}(\lambda):\lambda\in[\lambda_{1},\lambda_{2}]\}$ is VC-subgraph with dimension not greater than $2$.
However, as the centered VC-subgraph class is not necessarily VC-subgraph \citep[Problem 2.6.19]{VW23-MR4628026}, we need the following lemma to determine the VC dimension of the normalized class $\mathcal{F}$ whose proof can be found in \cref{pf:lem-VC-cov-residual}.

\begin{figure}
\centering
\includegraphics[width=\textwidth]{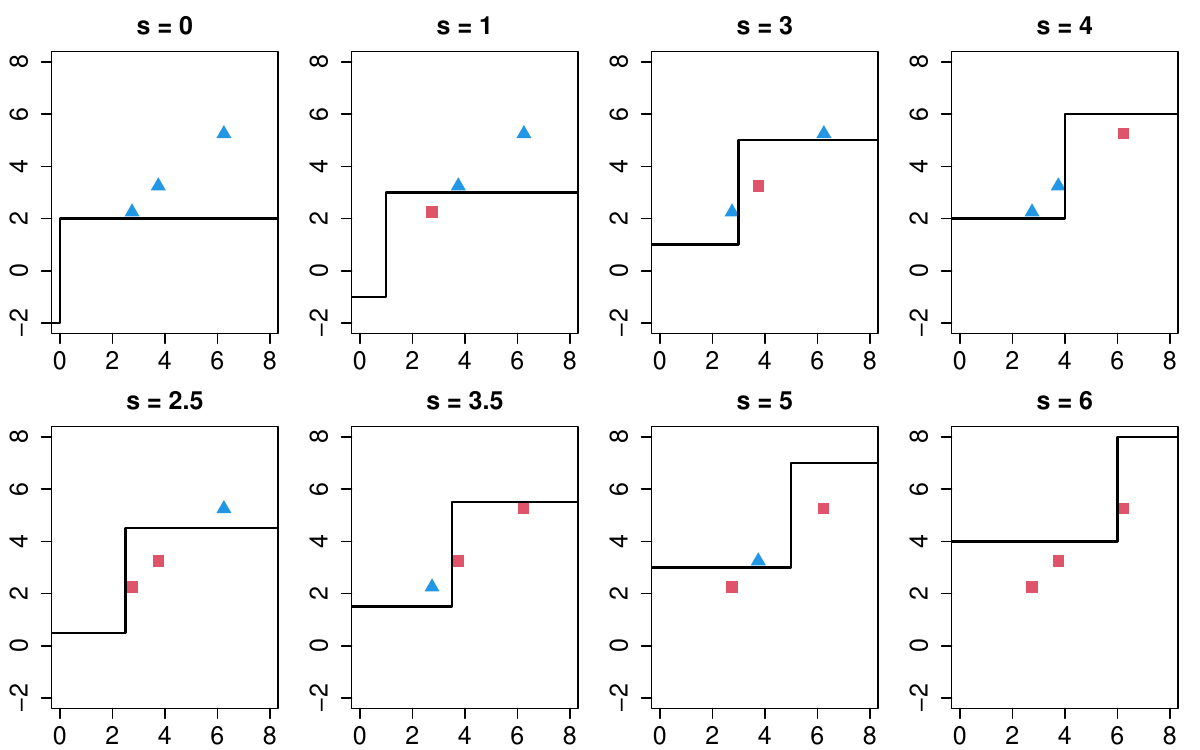}
\caption{Graphical illustration of $\vc(\mathcal{F})=3$ for the particular choice of $\mathcal{F}$ in \cref{lem:VC}(ii), given the three fixed points $(11/4,9/4)$, $(15/4,13/4)$, and $(25/4,21/4)$.
Each panel represents the selection of the subgraph of $f(\cdot,s)\in\mathcal{F}$ at a particular value of $s$, where the red squares indicate the points selected by the subgraph while the blue triangles are not selected.
}
\label{fig:VC}
\end{figure}

\begin{lemma}
\label{lem:VC}
(i) Let $F:\mathcal{X}\to\mathbb{R}$ be a nonnegative function.
Let $g:\mathcal{X}\times\mathbb{R}\to\mathbb{R}$ be a function such that $s\mapsto g(x,s)$ is nondecreasing for every $x\in\mathcal{X}$.
Let $h:\mathbb{R}\to\mathbb{R}$ be a nondecreasing function such that $h(s)\leq g(x,s)$ for every $x\in\mathcal{X}$ and $s\in\mathbb{R}$.
Let $\mathcal{F}=\{f(\cdot,s):s\in\mathbb{R}\}$ be a function class, where
\begin{equation}
\label{eq:lem-VC}
f(x,s)=g(x,s)1_{[s,\infty)}(F(x))+h(s)1_{[0,s)}(F(x)).
\end{equation}
Then $\mathcal{F}$ is a VC subgraph class with $\vc(\mathcal{F})\leq3$.

(ii) In particular, if we take $\mathcal{X}=[0,\infty)$, $F(x)=x$, $g(x,s)=s+2$, and $h(s)=s-2$, then $\vc(\mathcal{F})=3$ as illustrated in \cref{fig:VC}.
\end{lemma}

%
%

\iffalse
\begin{remark}
\cref{lem:VC}(i) is a nontrivial generalization of Lemma 2.6.18 and Example 2.6.24 in \cite{VW23-MR4628026}.
Intuitively, the component of $x\mapsto g(x,s)1_{[s,\infty)}(F(x))$ in \eqref{eq:lem-VC} is generalized from those considered in the mentioned textbook, which are VC of dimensions no more than 2.
The additional component of $x\mapsto h(s)1_{[0,s)}(F(x))$ over $s\in\mathbb{R}$ is also VC of dimension $1$ by an argument similar to Lemma 2.6.17 in \cite{VW23-MR4628026}.
Then by ``adding up'' these two components, we may have $\vc(\mathcal{F})\leq2+1=3$.
See \cref{pf:lem-VC} for a rigorous proof.
Furthermore, \cref{lem:VC}(ii) demonstrates that the upper bound of $\vc(\mathcal{F})\leq3$ is sharp.
\end{remark}
\fi
%
%

Indeed, observe that we may rewrite $f_{\lambda}$ as
\begin{equation}
\label{eq:HC-VC-decomposition}
f_{\lambda}(x)=\frac{\sigma_{0}(\lambda)}{\pi_{0}(\lambda)}1_{[\lambda,\infty)}(|x|)-\frac{\pi_{0}(\lambda)}{\sigma_{0}(\lambda)}1_{[0,\lambda)}(|x|).
\end{equation}
As $\lambda\mapsto\sigma_{0}(\lambda)/\pi_{0}(\lambda)$ is increasing over $[\lambda_{1},\lambda_{2}]$, it follows from \cref{lem:VC} that $\vc(\mathcal{F})\leq3$ for the higher criticism function class.

\subsubsection*{Step 2: Determining $\delta_{k}$ for $k=2,3,4$}
To determine $\delta_{2}$, $\delta_{3}$, and $\delta_{4}$ in \cref{thm:GA-mixing}, we note that by \eqref{eq:HC-VC-decomposition}, for $\lambda\in[\lambda_{1},\lambda_{2}]$ and $k=2,3,4$,
\begin{equation}
\Phi|f_{\lambda}|^{k}=\frac{\sigma_{0}^{k}(\lambda)}{\pi_{0}^{k}(\lambda)}\pi_{0}(\lambda)+\frac{\pi_{0}^{k}(\lambda)}{\sigma_{0}^{k}(\lambda)}(1-\pi_{0}(\lambda))\leq\frac{2\sigma_{0}^{k}(\lambda)}{\pi_{0}^{k-1}(\lambda)}\leq\frac{2}{\sigma_{0}^{k-2}(\lambda)}\leq2b^{k-2}.
\end{equation}
Then we may take $\delta_{k}^{k}=2/b^{2}$ for $k=2,3,4$ and thus $b\delta_{k}=2^{1/k}b^{1-2/k}$.

\subsubsection*{Step 3: Determining $\|\sigma_{P_{0}}-\sigma_{P_{1}}\|_{\mathcal{F}^{2}}$}
It remains to estimate $\|\sigma_{P_{0}}-\sigma_{P_{1}}\|_{\mathcal{F}^{2}}$.
Let $\Phi(x,y;\rho)=\mathbb{P}(X\leq x,Y\leq y)$ and $\phi(x,y;\rho)$ be, respectively, the distribution and density of $(X,Y)\sim N(0,(\begin{smallmatrix}1 & \rho\\
\rho & 1
\end{smallmatrix}))$.
Let $\Psi(x,y;\rho)=\mathbb{P}(|X|\geq x,|Y|\geq y)=2\Phi(-x,-y;\rho)+2\Phi(-x,-y;-\rho)$ and $\psi(x,y;\rho)=2\phi(x,y;\rho)+2\phi(x,y;-\rho)$.
Then we have for $\lambda_{1}\leq\lambda\leq\nu\leq\lambda_{2}$,
\begin{align}
&\sigma_{0}(\lambda)\sigma_{0}(\nu)\sigma_{P_{0}}(f_{\lambda},f_{\nu})=\sigma_{0}(\lambda)\sigma_{0}(\nu)\mathbb{E}(\mathbb{G}_{p}f_{\lambda}\mathbb{G}_{p}f_{\nu})\notag\\&=\frac{1}{p}\sum_{i=1}^{p}\sum_{j=1}^{p}(\mathbb{E}1_{[\lambda,\infty)}(|T_{i}|)1_{[\nu,\infty)}(|T_{j}|)-\pi_{0}(\lambda)\pi_{0}(\nu))\notag\\&=\pi_{0}(\nu)(1-\pi_{0}(\lambda))+2\sum_{k=1}^{p-1}\frac{p-k}{p}(\Psi(\lambda,\nu;\rho_{k})-\Psi(\lambda,\nu;0))\notag\\&=\pi_{0}(\nu)(1-\pi_{0}(\lambda))+2\sum_{k=1}^{p-1}\frac{p-k}{p}\rho_{k}\psi(\lambda,\nu;\theta_{k}\rho_{k})\notag\\&=\pi_{0}(\nu)(1-\pi_{0}(\lambda))+O(1)\sup_{|\rho|\leq\rho_{0}}\phi(\lambda,\nu;\rho),
\end{align}
where $\theta_{k}\in(0,1)$ and $|O(1)|\leq8\sum_{k=1}^{\infty}|\rho_{k}|$.
Here, we have used the following lemma whose proof can be found in \cref{pf:lem-VC-cov-residual}.

\begin{lemma}
\label{lem:cdf-rho}
We have $\partial_{\rho}\Phi(x,y;\rho)=\phi(x,y;\rho)$ and thus $\partial_{\rho}\Psi(x,y;\rho)=\psi(x,y;\rho)$.
\end{lemma}

Consequently, we have
\begin{equation}
\label{eq:HC-cov}
\sigma_{P_{0}}(f_{\lambda},f_{\nu})=\frac{\pi_{0}(\nu)(1-\pi_{0}(\lambda))}{\sigma_{0}(\lambda)\sigma_{0}(\nu)}+O(1)\sup_{|\rho|\leq\rho_{0}}\frac{\phi(\lambda,\nu;\rho)}{\sigma_{0}(\lambda)\sigma_{0}(\nu)},
\end{equation}
and similarly,
\begin{equation}
\sigma_{P_{1}}(f_{\lambda},f_{\nu})=\frac{ma_{1}}{p}\frac{\pi_{0}(\nu)(1-\pi_{0}(\lambda))}{\sigma_{0}(\lambda)\sigma_{0}(\nu)}+O(1)\sup_{|\rho|\leq\rho_{0}}\frac{\phi(\lambda,\nu;\rho)}{\sigma_{0}(\lambda)\sigma_{0}(\nu)}.
\end{equation}
By the Cauchy--Schwarz inequality, $\pi_{0}(\nu)(1-\pi_{0}(\lambda))\leq\sigma_{0}(\lambda)\sigma_{0}(\nu)$.
Then by the next lemma and $p\leq(m+1)(a_{1}+a_{2})$, we can conclude that
\begin{equation}
\label{eq:cov-function-approx}
\|\sigma_{P_{0}}-\sigma_{P_{1}}\|_{\mathcal{F}^{2}}=O\Big(\frac{1}{m}+\frac{a_{2}}{a_{1}}+\frac{(\log\log p)^{(1-\rho_{0})/2(1+\rho_{0})}}{(\log p)^{(1-\rho_{0})d/(1+\rho_{0})-1/2}}\Big).
\end{equation}

\begin{lemma}
\label{lem:cov-residual}
We have
\begin{equation}
\sup_{\lambda,\nu\in[\lambda_{1},\lambda_{2}]}\sup_{|\rho|\leq\rho_{0}}\frac{\phi(\lambda,\nu;\rho)}{\sigma_{0}(\lambda)\sigma_{0}(\nu)}=O\Big(\frac{(\log\log p)^{(1-\rho_{0})/2(1+\rho_{0})}}{(\log p)^{(1-\rho_{0})d/(1+\rho_{0})-1/2}}\Big).
\end{equation}
\end{lemma}

The proof of \cref{lem:cov-residual} can also be found in \cref{pf:lem-VC-cov-residual}.
Note that \eqref{eq:HC-almost-independent} in \cref{thm:HC-limit} is a direct consequence of \eqref{eq:HC-cov} and \cref{lem:cov-residual}.

\subsubsection*{Step 4: Conclusion}
Let $a_{1}=(\log p)^{\tau_{1}}$ and $a_{2}=(\log p)^{\tau_{2}}$ for $\tau_{1}\geq\tau_{2}>1/c_{2}$.
Then we can apply \cref{thm:GA-mixing} with \cref{rmk:thm-GA-mixing} to conclude that there is a random variable $\HC_{p}^{**}$ as described in \cref{thm:GA-HC} such that
\begin{equation}
\mathbb{P}(|\HC_{p}^{*}-\HC_{p}^{**}|>\Delta_{0})\leq C_{q}\gamma+o(1),
\end{equation}
where
\begin{align}
\gamma\Delta_{0}&=O((\log p)^{(2+2\tau_{1}-c)/2}+(\log p)^{(3+3\tau_{1}-c)/4}+(\log p)^{(4+4\tau_{1}-c)/6})\notag\\&\quad+O(p^{-1/2}(\log p)^{(1+\tau_{1})/2}+(\log p)^{(1-\tau_{1}+\tau_{2})/2})\notag\\&\quad+O\Big(\frac{(\log\log p)^{(1-\rho_{0})/4(1+\rho_{0})}}{(\log p)^{(1-\rho_{0})d/2(1+\rho_{0})-3/4}}+(\log p)^{(1-\tau_{1}+2\tau_{2})/2}\Big)\notag\\&=O((\log p)^{-\tau}).
\label{eq:Delta-0-HC}
\end{align}
Here, we may take
\begin{equation}
0<\tau<\frac{c-4\tau_{1}-4}{6}\wedge\frac{\tau_{1}-2\tau_{2}-1}{2}\wedge\frac{2(1-\rho_{0})d-3(1+\rho_{0})}{4(1+\rho_{0})},
\end{equation}
which implies \eqref{eq:thm-GA-HC} in \cref{thm:GA-HC}.

To prove the remaining part of \cref{thm:HC-limit}, let $\sigma_{P_{0}}^{0}(f_{\lambda},f_{\nu})=\pi_{0}(\nu)(1-\pi_{0}(\lambda))/\sigma_{0}(\lambda)\sigma_{0}(\nu)$ be the covariance function of $\mathbb{B}_{p}^{0}$ for $\lambda_{1}\leq\lambda\leq\nu\leq\lambda_{2}$.
Then \eqref{eq:HC-almost-independent} in \cref{thm:HC-limit} can be rewritten as
\begin{equation}
\label{eq:cov-function-approx-BB}
\|\sigma_{P_{0}}-\sigma_{P_{0}}^{0}\|_{\mathcal{F}^{2}}=O\Big(\frac{(\log\log p)^{(1-\rho_{0})/2(1+\rho_{0})}}{(\log p)^{(1-\rho_{0})d/(1+\rho_{0})-1/2}}\Big).
\end{equation}
By \cref{lem:comparison-process} with a slight modification in the final step of the proof of \cref{thm:GA-mixing}, there is a random variable $\BB_{p}^{*}$ as described in \cref{sec:intro} such that
\begin{equation}
\mathbb{P}(|\HC_{p}^{*}-\BB_{p}^{*}|>\Delta_{0})\leq C_{q}\gamma+o(1),
\end{equation}
where $\Delta_{0}$ has been defined in \eqref{eq:Delta-0-HC}.
This implies \eqref{eq:HC-significance-level-GA}.
Then \eqref{eq:HC-significance-level-asymptotics} is a direct consequence of Lemma 1 in \cite{Jae79-MR0515687}.
This completes the proof.

\section{Proof of \cref{thm:GA-HC-t}}
\label{pf:thm-GA-HC-t}
The proof is a modification of the previous proof.
Note that $\HC_{p}^{*}$ itself is not ready the supremum of an empirical process in this case due to the non-Gaussianity, which requires an additional approximation.

\subsubsection*{Step 1: Approximating $\HC_{p}^{*}$ by the supremum of an empirical process}
Observe that the higher criticism $t$-statistic can be rewritten as
\begin{equation}
\label{eq:HC-threhsolding-t}
\HC_{p}^{*}=\sup_{\lambda\in[\lambda_{1},\lambda_{2}]}\frac{1}{p^{1/2}\sigma_{0}(\lambda)}\sum_{j=1}^{p}(1_{[\lambda,\infty)}(|T_{j}|)-\pi_{0}(\lambda)).
\end{equation}
Let $P$ be the distribution of $T_{j}$.
Let $\mathbb{G}_{p}=p^{-1/2}\sum_{j=1}^{p}(\delta_{T_{j}}-P)$ be the empirical process associated to $T_{1},\dots,T_{p}$.
Let $\pi_{1}(\lambda)=\mathbb{P}(|T_{j}|\geq\lambda)$ for $\lambda\geq0$.
Then \eqref{eq:HC-threhsolding-t} can be approximated by
\begin{equation}
\HC_{p}^{**}=\sup_{\lambda\in[\lambda_{1},\lambda_{2}]}\mathbb{G}_{p}f_{\lambda,1},\quad\text{where}\quad f_{\lambda,1}(x)=\frac{1}{\sigma_{0}(\lambda)}(1_{[\lambda,\infty)}(|x|)-\pi_{1}(\lambda)).
\end{equation}
Indeed, we observe that
\begin{equation}
\mathbb{E}_{p}f_{\lambda}(T_{j})=\mathbb{E}_{p}f_{\lambda,1}(T_{j})+\frac{\pi_{1}(\lambda)-\pi_{0}(\lambda)}{\sigma_{0}(\lambda)},
\end{equation}
where $f_{\lambda}$ has been defined in \eqref{eq:HC-empirical-process} and $\mathbb{E}_{p}a_{j}=p^{-1}\sum_{j=1}^{p}a_{j}$.
By \eqref{eq:MDR-1} in \cref{lem:MDR},
\begin{equation}
\frac{\pi_{1}(\lambda)-\pi_{0}(\lambda)}{\sigma_{0}(\lambda)}=\pi_{0}^{1/2}(\lambda)O\Big(\frac{(\log p)^{1+\delta/2}}{n^{\delta/2}}\Big)=O\Big(\frac{(\log p)^{1+\delta/2}}{p^{d/2}n^{\delta/2}}\Big).
\end{equation}
uniformly in $\lambda\in[\lambda_{1},\lambda_{2}]$.
This implies the following estimate of the approximation error bound and it suffices to investigate $\HC_{p}^{**}$.
\begin{equation}
\label{eq:HC-t-approximation}
|\HC_{p}^{*}-\HC_{p}^{**}|=O\Big(\frac{p^{(1-d)/2}(\log p)^{1+\delta/2}}{n^{\delta/2}}\Big)
\end{equation}

\subsubsection*{Step 2: Approximating $\HC_{p}^{**}$ by $\BB_{p}^{**}$}
Let $\sigma_{1}^{2}(\lambda)=\pi_{1}(\lambda)(1-\pi_{1}(\lambda))$.
Let $\mathcal{F}_{1}=\{f_{\lambda,1}:\lambda\in[\lambda_{1},\lambda_{2}]\}$ and $F_{1}(x)=1/(\sigma_{0}(\lambda_{2})\wedge\sigma_{1}(\lambda_{2}))$.
Let $\mathbb{B}_{p}$ be a tight Gaussian process in $\ell^{\infty}(\mathcal{F}_{1})$ with mean zero and covariance function $\mathbb{E}(\mathbb{B}_{p}f_{\lambda,1}\mathbb{B}_{p}f_{\nu,1})=\mathbb{E}(\mathbb{G}_{p}f_{\lambda,1}\mathbb{G}_{p}f_{\nu,1})$.
Then we shall approximate $\mathbb{G}_{p}$ by $\mathbb{B}_{p}$.
To apply \cref{thm:GA-mixing}, by the proof in \cref{pf:thm-GA-HC} we similarly have $b=1/(\sigma_{0}(\lambda_{2})\wedge\sigma_{1}(\lambda_{2}))$, $q=\infty$, $v=3$, and $b\delta_{k}=2^{1/k}b^{1-2/k}$ for $k=2,3,4$.
It remains to determine $\|\sigma_{P_{0}}-\sigma_{P_{1}}\|_{\mathcal{F}_{1}^{2}}$.
To this end, we shall use the following inequalities that can be found as Proposition 1.1.1 and Lemma 1.2.3 in \cite{Dou94-MR1312160}.
Let $\alpha(\mathcal{A},\mathcal{B})=\sup\{|\mathbb{P}(A)\mathbb{P}(B)-\mathbb{P}(A\cap B)|:A\in\mathcal{A},B\in\mathcal{B}\}$ be the $\alpha$-mixing coefficient between $\sigma$-algebras $\mathcal{A}$ and $\mathcal{B}$.

\begin{lemma}
\label{lem:cov-inequality}
(i) We have $2\alpha(\mathcal{A},\mathcal{B})\leq\beta(\mathcal{A},\mathcal{B})$.
(ii) Let $X$ and $Y$ be two random variables.
Then we have $|\cov(X,Y)|\leq4\alpha(\sigma(X),\sigma(Y))\|X\|_{\infty}\|Y\|_{\infty}$.
\end{lemma}

By Lemmas \ref{lem:MDR} and \ref{lem:cov-inequality} and a direct calculation, we have for $\lambda_{1}\leq\lambda\leq\nu\leq\lambda_{2}$,
\begin{align}
&\sigma_{0}(\lambda)\sigma_{0}(\nu)\sigma_{P_{0}}(f_{\lambda,1},f_{\nu,1})=\sigma_{0}(\lambda)\sigma_{0}(\nu)\mathbb{E}(\mathbb{G}_{p}f_{\lambda,1}\mathbb{G}_{p}f_{\nu,1})\notag\\&=\frac{1}{p}\sum_{j=1}^{p}\sum_{k=1}^{p}\cov(1_{[\lambda,\infty)}(|T_{j}|),1_{[\nu,\infty)}(|T_{k}|))\notag\\&=\cov(1_{[\lambda,\infty)}(|T_{1}|),1_{[\nu,\infty)}(|T_{1}|))\notag\\&\quad+2\Big(\sum_{k=1}^{h-1}+\sum_{k=h}^{p-1}\Big)\frac{p-k}{p}\cov(1_{[\lambda,\infty)}(|T_{1}|),1_{[\nu,\infty)}(|T_{1+k}|))\notag\\&=\mathbb{E}1_{[\lambda,\infty)}(|T_{1}|)1_{[\nu,\infty)}(|T_{1}|)-\pi_{1}(\lambda)\pi_{1}(\nu)+O(p\beta(h))\notag\\&\quad+2\sum_{k=1}^{h-1}\frac{p-k}{p}(\mathbb{E}1_{[\lambda,\infty)}(|T_{1}|)1_{[\nu,\infty)}(|T_{1+k}|)-\pi_{1}(\lambda)\pi_{1}(\nu))\notag\\&=\pi_{0}(\nu)(1-\pi_{0}(\lambda))\Big(1+O\Big(\frac{(\log p)^{1+\delta/2}}{n^{\delta/2}}\Big)\Big)+\pi_{0}(\lambda)\pi_{0}(\nu)O\Big(\frac{h(\log p)^{1+\delta/2}}{n^{\delta/2}}\Big)\notag\\&\quad+2\sum_{k=1}^{h-1}\frac{p-k}{p}(\Psi(\lambda,\nu;\rho_{k})-\Psi(\lambda,\nu;0))\Big(1+O\Big(\frac{(\log p)^{1+\delta/2}}{n^{\delta/2}}\Big)\Big)+O(p\beta(h))\notag\\&=\pi_{0}(\nu)(1-\pi_{0}(\lambda))\notag\\&\quad+O\Big(\sup_{|\rho|\leq\rho_{0}}\phi(\lambda,\nu;\rho)+p\beta(h)+\frac{\pi_{0}(\lambda)\pi_{0}(\nu)h(\log p)^{1+\delta/2}}{n^{\delta/2}}\Big).\label{eq:cov-function-approx-t-1}
\end{align}
To estimate the term $\sup\{\phi(\lambda,\nu;\rho):|\rho|\leq\rho_{0}\}$ in this case, we shall apply the following lemma whose proof is similar to that of \cref{lem:cov-residual} and is omitted.

\begin{lemma}
\label{lem:cov-residual-t}
We have
\begin{equation}
\sup_{\lambda,\nu\in[\lambda_{1},\lambda_{2}]}\sup_{|\rho|\leq\rho_{0}}\frac{\phi(\lambda,\nu;\rho)}{\sigma_{0}(\lambda)\sigma_{0}(\nu)}=O\Big(\frac{(\log p)^{1/(1+\rho_{0})}}{p^{(1-\rho_{0})d/(1+\rho_{0})}}\Big).
\end{equation}
\end{lemma}

By a similar argument as that in deriving the estimate in \eqref{eq:cov-function-approx} and \eqref{eq:cov-function-approx-BB}, we have for $\tau_{1}\geq\tau_{2}>1/c_{2}$ and $h=((3/c_{1})\log p)^{1/c_{2}}$,
\begin{align}
&\|\sigma_{P_{0}}-\sigma_{P_{1}}\|_{\mathcal{F}_{1}^{2}}\notag\\&=O\Big(\frac{(\log p)^{\tau_{1}}}{p}+\frac{1}{(\log p)^{\tau_{1}-\tau_{2}}}+\frac{(\log p)^{1/(1+\rho_{0})}}{p^{(1-\rho_{0})d/(1+\rho_{0})}}+\frac{(\log p)^{1+\delta/2+1/c_{2}}}{p^{d}n^{\delta/2}}\Big),\label{eq:cov-function-approx-t-2}
\end{align}
and
\begin{equation}
\|\sigma_{P_{0}}-\sigma_{P_{0}}^{0}\|_{\mathcal{F}_{1}^{2}}=O\Big(\frac{(\log p)^{1/(1+\rho_{0})}}{p^{(1-\rho_{0})d/(1+\rho_{0})}}+\frac{(\log p)^{1+\delta/2+1/c_{2}}}{p^{d}n^{\delta/2}}\Big),
\end{equation}
where $\sigma_{P_{0}}^{0}(f_{\lambda,1},f_{\nu,1})=\pi_{0}(\nu)(1-\pi_{0}(\lambda))/\sigma_{0}(\lambda)\sigma_{0}(\nu)$.
Consequently, we can apply \cref{thm:GA-mixing} with \cref{rmk:thm-GA-mixing} and \cref{lem:comparison-process} to obtain that there is a random variable $\BB_{p}^{*}$ described as in \cref{sec:intro} such that
\begin{equation}
\mathbb{P}(|\HC_{p}^{**}-\BB_{p}^{*}|>\Delta_{0})\leq C_{q}\gamma+o(1),
\end{equation}
where
\begin{align}
\gamma\Delta_{0}&=O((\log p)^{(2+2\tau_{1}-c)/2}+(\log p)^{(3+3\tau_{1}-c)/4}+(\log p)^{(4+4\tau_{1}-c)/6})\notag\\&\quad+O\Big(\frac{(\log p)^{(1+\tau_{1})/2}}{p^{1/2}}+(\log p)^{(1-\tau_{1}+\tau_{2})/2}+\frac{(\log p)^{(2+\rho_{0})/2(1+\rho_{0})}}{p^{(1-\rho_{0})d/2(1+\rho_{0})}}\Big)\notag\\&\quad+O\Big(\frac{(\log p)^{1+\delta/4+1/2c_{2}}}{p^{d/2}n^{\delta/4}}+(\log p)^{(1-\tau_{1}+2\tau_{2})/2}\Big)\notag\\&=O((\log p)^{-\tau}),
\end{align}
and $0<\tau<((c-4\tau_{1}-4)/6)\wedge((\tau_{1}-2\tau_{2}-1)/2)$.
This implies that
\begin{equation}
\label{eq:HC-t-approximation-2}
|\HC_{p}^{**}-\BB_{p}^{*}|=O_{\mathbb{P}}((\log p)^{-\tau}).
\end{equation}
Then \eqref{eq:thm-GA-HC-t-1} follows from \eqref{eq:HC-t-approximation} and \eqref{eq:HC-t-approximation-2}.
The second statement of \cref{thm:GA-HC-t} follows from a direct application of Lemma 1 in \cite{Jae79-MR0515687}.

\section{Proof of \cref{thm:GA-MT-t}}
\label{pf:thm-GA-MT-t}
The proof is similar to the previous proof with the help of the moment calculation results in \cite{QCS25+} for each $\MT_{p,\lambda}$ and the limit result \cref{lem:Hus90-modification} for the supremum of $\mathbb{B}_{p}$.

\subsubsection*{Step 1: Approximating $\MT_{p}^{*}$ by the supremum of an empirical process}
Recall that $\mathbb{G}_{p}=p^{-1/2}\sum_{j=1}^{p}(\delta_{T_{j}}-P)$ is the empirical process associated to $T_{1},\dots,T_{p}$.
Let $\mu_{1}(\lambda)=\mathbb{E}|T_{j}|^{2}1_{[\lambda,M]}(|T_{j}|)$ and $\sigma_{1}^{2}(\lambda)=\mathbb{E}|T_{j}|^{4}1_{[\lambda,M]}(|T_{j}|)-\mu_{1}^{2}(\lambda)$ for $\lambda\geq0$.
Then \eqref{eq:MT-threhsolding-t} can be approximated by
\begin{equation}
\MT_{p}^{**}=\sup_{\lambda\in[\lambda_{1},\lambda_{2}]}\mathbb{G}_{p}f_{\lambda,1},\ \text{where}\ f_{\lambda,1}(x)=\frac{1}{\sigma_{0}(\lambda)}(|x|^{2}1_{[\lambda,M]}(|x|)-\mu_{1}(\lambda)).
\end{equation}
By \eqref{eq:MDR-1} in \cref{lem:MDR} and a similar argument as that in proving Theorem 3.1 of \cite{QCS25+},
\begin{equation}
\frac{\mu_{1}(\lambda)-\mu_{0}(\lambda)}{\sigma_{0}(\lambda)}=\frac{\mu_{0}(\lambda)}{\sigma_{0}(\lambda)}O\Big(\frac{(\log p)^{1+\delta/2}}{n^{\delta/2}}\Big)=O\Big(\frac{(\log p)^{\delta/2}}{p^{d/2}n^{\delta/2}}\Big).
\end{equation}
uniformly in $\lambda\in[\lambda_{1},\lambda_{2}]$.
Then we have
\begin{equation}
\label{eq:thm-GA-MT-1}
|\MT_{p}^{*}-\MT_{p}^{**}|=O\Big(\frac{p^{(1-d)/2}(\log p)^{\delta/2}}{n^{\delta/2}}\Big).
\end{equation}

\subsubsection*{Step 2: Approximating $\MT_{p}^{**}$ by $\mathbb{B}_{p}^{*}$}
Let $\mathcal{F}_{1}=\{f_{\lambda,1}:\lambda\in[\lambda_{1},\lambda_{2}]\}$ and $F_{1}(x)=2M^{2}/\sigma_{0}(\lambda_{2})$.
Let $\mathbb{B}_{p}$ be a tight Gaussian process in $\ell^{\infty}(\mathcal{F}_{1})$ with mean zero and covariance function given by $\mathbb{E}(\mathbb{B}_{p}f_{\lambda,1}\mathbb{B}_{p}f_{\nu,1})=\mathbb{E}(\mathbb{G}_{p}f_{\lambda,1}\mathbb{G}_{p}f_{\nu,1})$.
Then we shall approximate $\mathbb{G}_{p}$ by $\mathbb{B}_{p}$.
To apply \cref{thm:GA-mixing}, by the proof in Sections \ref{pf:thm-GA-HC} and \ref{pf:thm-GA-HC-t} we similarly have $b=2M^{2}/\sigma_{0}(\lambda_{2})$, $q=\infty$, $v=3$, and $b\delta_{k}=b^{1-2/k}$ for $k=2,3,4$ and $p$ sufficiently large.
It remains to determine $\|\sigma_{P_{0}}-\sigma_{P_{1}}\|_{\mathcal{F}_{1}^{2}}$.
Indeed, for $\lambda_{1}\leq\lambda\leq\nu\leq\lambda_{2}$,
\begin{align}
&\sigma_{0}(\lambda)\sigma_{0}(\nu)\sigma_{P_{0}}(f_{\lambda,1},f_{\nu,1})=\sigma_{0}(\lambda)\sigma_{0}(\nu)\mathbb{E}(\mathbb{G}_{p}f_{\lambda,1}\mathbb{G}_{p}f_{\nu,1})\notag\\&=\frac{1}{p}\sum_{j=1}^{p}\sum_{k=1}^{p}\cov(|T_{j}|^{2}1_{[\lambda,M]}(|T_{j}|),|T_{k}|^{2}1_{[\nu,M]}(|T_{k}|))\notag\\&=\cov(|T_{1}|^{2}1_{[\lambda,M]}(|T_{1}|),|T_{1}|^{2}1_{[\nu,M]}(|T_{1}|))\notag\\&\quad+2\Big(\sum_{k=1}^{h-1}+\sum_{k=h}^{p-1}\Big)\frac{p-k}{p}\cov(|T_{1}|^{2}1_{[\lambda,M]}(|T_{1}|),|T_{1+k}|^{2}1_{[\nu,M]}(|T_{1+k}|))\notag\\&=\mathbb{E}|T_{1}|^{4}1_{[\lambda,M]}(|T_{1}|)1_{[\nu,M]}(|T_{1}|)-\mu_{1}(\lambda)\mu_{1}(\nu)+O(M^{4}p\beta(h))\notag\\&\quad+2\sum_{k=1}^{h-1}\frac{p-k}{p}(\mathbb{E}|T_{1}|^{2}|T_{1+k}|^{2}1_{[\lambda,M]}(|T_{1}|)1_{[\nu,M]}(|T_{1+k}|)-\mu_{1}(\lambda)\mu_{1}(\nu))\notag\\&=(\sigma_{0}^{2}(\nu)+\mu_{0}(\nu)(\mu_{0}(\nu)-\mu_{0}(\lambda)))\Big(1+O\Big(\frac{(\log p)^{1+\delta/2}}{n^{\delta/2}}\Big)\Big)\notag\\&\quad+2\sum_{k=1}^{h-1}\frac{p-k}{p}(\Xi(\lambda,\nu;\rho_{k})-\Xi(\lambda,\nu;0))\Big(1+O\Big(\frac{(\log p)^{1+\delta/2}}{n^{\delta/2}}\Big)\Big)\notag\\&\quad+\mu_{0}(\lambda)\mu_{0}(\nu)O\Big(\frac{h(\log p)^{1+\delta/2}}{n^{\delta/2}}\Big)+O(M^{4}p\beta(h))\notag\\&=(\sigma_{0}^{2}(\nu)+\mu_{0}(\nu)(\mu_{0}(\nu)-\mu_{0}(\lambda)))\notag\\&\quad+O\Big(\frac{L_{p}\sigma_{0}(\lambda)\sigma_{0}(\nu)}{p^{(1-\rho_{0})d/(1+\rho_{0})}}+M^{4}p\beta(h)+\frac{M^{4}h(\log p)^{1+\delta/2}}{n^{\delta/2}}\Big).
\label{eq:cov-MT}
\end{align}
Here, for $(Z_{1},Z_{1+k})\sim N(0,(\begin{smallmatrix}1 & \rho_{k}\\
\rho_{k} & 1
\end{smallmatrix}))$, by Fubini's theorem we have
\begin{align}
&\Xi(\lambda,\nu;\rho_{k})=\mathbb{E}|Z_{1}|^{2}|Z_{1+k}|^{2}1_{[\lambda,M]}(|Z_{1}|)1_{[\nu,M]}(|Z_{1+k}|)\notag\\&=\mathbb{E}\int_{0}^{|Z_{1}|^{2}1_{[\lambda,M]}(|Z_{1}|)}dx\int_{0}^{|Z_{1+k}|^{2}1_{[\nu,M]}(|Z_{1+k}|)}dy\notag\\&=\int_{0}^{\infty}\int_{0}^{\infty}\mathbb{P}(|Z_{1}|^{2}1_{[\lambda,M]}(|Z_{1}|)\geq x,|Z_{1+k}|^{2}1_{[\nu,M]}(|Z_{1+k}|)\geq y)dxdy\notag\\&=\lambda^{2}\nu^{2}\mathbb{P}(|Z_{1}|\in[\lambda,M],|Z_{1+k}|\in[\nu,M])\notag\\&\quad+\nu^{2}\int_{0}^{\infty}\mathbb{P}(|Z_{1}|\in[x^{1/2},M],|Z_{1+k}|\in[\nu,M])dx\notag\\&\quad+\lambda^{2}\int_{0}^{\infty}\mathbb{P}(|Z_{1}|\in[\lambda,M],|Z_{1+k}|\in[y^{1/2},M])dy\notag\\&\quad+\int_{0}^{\infty}\int_{0}^{\infty}\mathbb{P}(|Z_{1}|\in[x^{1/2},M],|Z_{1+k}|\in[y^{1/2},M])dxdy.
\end{align}
Each probability in the last equation can be written as a four-term summation of the $\Psi$ function defined in \cref{pf:thm-GA-HC}.
Then by successively applying the mean value theorem to each pair of $\Psi(\cdot,\cdot;\rho_{k})-\Psi(\cdot,\cdot;0)$ and a similar argument as that in \cref{lem:cov-residual-t}, we can conclude that
\begin{equation}
\sup_{\lambda,\nu\in[\lambda_{1},\lambda_{2}]}\frac{|\Xi(\lambda,\nu;\rho_{k})-\Xi(\lambda,\nu;0)|}{\sigma_{0}(\lambda)\sigma_{0}(\nu)}=O\Big(\frac{|\rho_{k}|L_{p}}{p^{(1-\rho_{0})d/(1+\rho_{0})}}\Big),
\end{equation}
where $L_{p}$ is a polylogarithmic term.
Let $\sigma_{P_{0}}^{0}(f_{\lambda,1},f_{\nu,1})=\mathbb{E}(\mathbb{B}_{p,\lambda}\mathbb{B}_{p,\nu})$.
Then by \eqref{eq:cov-brownian-bridge-MT}, \eqref{eq:cov-MT}, and a similar argument as that in \cref{pf:thm-GA-HC-t}, we have
\begin{equation}
\|\sigma_{P_{0}}-\sigma_{P_{0}}^{0}\|_{\mathcal{F}_{1}^{2}}=O\Big(\frac{L_{p}}{p^{(1-\rho_{0})d/(1+\rho_{0})}\wedge(p^{d}n^{\delta/2})}\Big),
\end{equation}
and
\begin{equation}
\label{eq:thm-GA-MT-2}
|\MT_{p}^{**}-\mathbb{B}_{p}^{*}|=O((\log p)^{-\tau}),
\end{equation}
where $\tau$ is the same as that in \eqref{eq:HC-t-approximation-2}.
Then \eqref{eq:thm-GA-MT} follows from \eqref{eq:thm-GA-MT-1} and \eqref{eq:thm-GA-MT-2}.

\subsubsection*{Step 3: Deriving asymptotic distribution of $\mathbb{B}_{p}^{*}$}
Observe that $\lambda\mapsto\sigma_{0}^{2}(\lambda)$ is strictly decreasing in $[\lambda_{1},\lambda_{2}]$ for $p$ sufficiently large.
Then we may represent $\mathbb{B}_{p}^{*}$ as
\begin{equation}
\label{eq:multi-level-thresholding-statistic-OU}
\mathbb{B}_{p}^{*}\overset{d}{=}\sup_{t\in[0,T_{p}]}\mathbb{X}_{p}(t),
\end{equation}
where $\mathbb{X}_{p}(t)=\mathbb{B}_{p,\lambda}$ for $t=\log\sigma_{0}^{2}(\lambda_{1})-\log\sigma_{0}^{2}(\lambda)$ and
\begin{equation}
T_{p}=\log\sigma_{0}^{2}(\lambda_{1})-\log\sigma_{0}^{2}(\lambda_{2})=(1-d)\log p-c\log\log p.
\end{equation}
In what follows, we shall use \cref{lem:Hus90-modification} to provide the asymptotic distribution of \eqref{eq:multi-level-thresholding-statistic-OU}.
Note that it suffices to show that $\{\mathbb{X}_{p}(t):t\in[0,T_{p}]\}$ is a sequence of asymptotically locally stationary Gaussian processes.

We shall investigate the covariance function of $\mathbb{B}_{p}$.
Note that
\begin{equation}
\mu_{0}(\lambda)\sim2\lambda\phi(\lambda)\quad\text{and}\quad\sigma_{0}^{2}(\lambda)\sim2\lambda^{3}\phi(\lambda)
\end{equation}
uniformly in $\lambda\in[\lambda_{p,1},\lambda_{p,2}]$.
This implies that
\begin{equation}
\frac{\mu_{0}(\nu)}{\sigma_{0}(\lambda)\sigma_{0}(\nu)}\sim\frac{\nu\phi(\nu)}{(\lambda^{3}\nu^{3}\phi(\lambda)\phi(\nu))^{1/2}}=\frac{1}{\nu^{2}}\Big(\frac{\nu^{3}\phi(\nu)}{\lambda^{3}\phi(\lambda)}\Big)^{1/2}=o\Big(\frac{\sigma_{0}(\nu)}{\sigma_{0}(\lambda)}\Big)
\label{eq:thm-MTT-iid-Gaussian-X-cov-2}
\end{equation}
uniformly in $\lambda_{1}\leq\lambda\leq\nu\leq\lambda_{2}$.
By combining \eqref{eq:cov-brownian-bridge-MT} and \eqref{eq:thm-MTT-iid-Gaussian-X-cov-2}, we have
\begin{align}
\mathbb{E}(\mathbb{B}_{p,\lambda}\mathbb{B}_{p,\nu})&=\exp\Big(-\frac{1}{2}|\log\sigma_{0}^{2}(\lambda)-\log\sigma_{0}^{2}(\nu)|\Big)\Big(1+o(|\mu_{0}(\lambda)-\mu_{0}(\nu)|)\Big)\notag\\&=1-\frac{1}{2}|\log\sigma_{0}^{2}(\lambda)-\log\sigma_{0}^{2}(\nu)|\notag\\&\quad+O(|\log\sigma_{0}^{2}(\lambda)-\log\sigma_{0}^{2}(\nu)|^{2})+o(|\mu_{0}(\lambda)-\mu_{0}(\nu)|)\label{eq:thm-MTT-iid-Gaussian-X-cov-3}
\end{align}
uniformly in $\lambda,\nu\in[\lambda_{1},\lambda_{2}]$.
By the mean value theorem,
\begin{align}
\frac{\mu_{0}(\lambda)-\mu_{0}(\nu)}{\log\sigma_{0}^{2}(\lambda)-\log\sigma_{0}^{2}(\nu)}&=\frac{\partial\mu_{0}(\lambda_{\theta})}{\partial\log\sigma_{0}^{2}(\lambda_{\theta})}\notag\\&=\frac{-2\lambda_{\theta}^{2}\phi(\lambda_{\theta})}{(4\lambda_{\theta}^{2}\phi(\lambda_{\theta})\mu_{0}(\lambda_{\theta})-2\lambda_{\theta}^{4}\phi(\lambda_{\theta}))/\sigma_{0}^{2}(\lambda_{\theta})}\notag\\&=\frac{\sigma_{0}^{2}(\lambda_{\theta})}{\lambda_{\theta}^{2}-2\mu_{0}(\lambda_{\theta})}\sim2\lambda_{\theta}\phi(\lambda_{\theta})\to0\label{eq:thm-MTT-iid-Gaussian-X-cov-4}
\end{align}
uniformly in $\lambda,\nu\in[\lambda_{1},\lambda_{2}]$, where $\lambda_{\theta}=\theta\lambda+(1-\theta)\nu$ for a value $\theta\in(0,1)$ that may depend on $\lambda$ and $\nu$.
By combining \eqref{eq:thm-MTT-iid-Gaussian-X-cov-3} and \eqref{eq:thm-MTT-iid-Gaussian-X-cov-4}, we can conclude that as $|\log\sigma_{0}^{2}(\lambda)-\log\sigma_{0}^{2}(\nu)|\to0$ and $p\to\infty$,
\begin{equation}
\mathbb{E}(\mathbb{B}_{p,\lambda}\mathbb{B}_{p,\nu})=1-\frac{1}{2}|\log\sigma_{0}^{2}(\lambda)-\log\sigma_{0}^{2}(\nu)|+o(|\log\sigma_{0}^{2}(\lambda)-\log\sigma_{0}^{2}(\nu)|)
\label{eq:thm-MTT-iid-Gaussian-X-cov-5}
\end{equation}
uniformly in $\lambda,\nu\in[\lambda_{1},\lambda_{2}]$,
and that as $\tau\to\infty$ and $p\to\infty$,
\begin{equation}
\mathbb{E}(\mathbb{B}_{p,\lambda}\mathbb{B}_{p,\nu})\log\tau\leq e^{-\tau/2}(1+o(\tau))\log\tau\to0
\label{eq:thm-MTT-iid-Gaussian-X-cov-6}
\end{equation}
uniformly in $\lambda,\nu\in[\lambda_{1},\lambda_{2}]$ and $|\log\sigma_{0}^{2}(\lambda)-\log\sigma_{0}^{2}(\nu)|\geq\tau$.

By the reparametrization in \eqref{eq:multi-level-thresholding-statistic-OU}, it can be deduced from \eqref{eq:thm-MTT-iid-Gaussian-X-cov-5} and \eqref{eq:thm-MTT-iid-Gaussian-X-cov-6} that as $|s-t|\to0$ and $p\to0$,
\begin{equation}
\mathbb{E}(\mathbb{X}_{p}(s)\mathbb{X}_{p}(t))=1-\frac{1}{2}|s-t|+o(|s-t|)
\end{equation}
uniformly in $s,t\in[0,T_{p}]$,
and that as $\tau\to\infty$ and $p\to\infty$,
\begin{equation}
\sup_{|s-t|\geq\tau}\mathbb{E}(\mathbb{X}_{p}(s)\mathbb{X}_{p}(t))\log\tau\to0.
\end{equation}
This shows that $\{\mathbb{X}_{p}(t):t\in[0,T_{p}]\}$ is a sequence of asymptotically locally stationary Gaussian processes satisfying \eqref{eq:asymptotically-locally-stationary-Gaussian-process} with parameters $\alpha=1$ and $C\equiv1/2$.
Furthermore, by \eqref{eq:cov-brownian-bridge-MT} and the monotonicity of $\mu_{0}$, we have $\mathbb{E}(\mathbb{X}_{p}(s)\mathbb{X}_{p}(t))=1$ if and only if $s=t$ for $p$ sufficiently large.
Then \eqref{eq:MT-limit} is a direct consequence of \eqref{eq:thm-GA-MT} and \cref{lem:Hus90-modification}.
This completes the proof.

\section{Proof of \cref{thm:CCK14}}
\label{pf:lem-CCK14}
We first introduce additional notations.
For a matrix $A=(a_{ij})$, let $\|A\|_{\max}=\max_{ij}|a_{ij}|$ be its maximum norm.
Let $\mathbb{E}_{n}a_{i}=n^{-1}\sum_{i=1}^{n}a_{i}$ for a sequence $a_{1},\dots,a_{n}$.
Let $\|X\|_{q}=(\mathbb{E}|X|^{q})^{1/q}$ for a real-valued random variable $X$.
Let $\mathcal{F}\cdot\mathcal{F}=\{fg:f\in\mathcal{F},g\in\mathcal{F}\}$.
Note that $\mathcal{F}\cdot\mathcal{F}$ should not be confused with $\mathcal{F}^{2}$.
Let $\|f\|_{\infty}=\sup_{x\in\mathcal{X}}|f(x)|$ be the supremum norm of $f$.
For $\delta\in(0,1)$, define the uniform entropy integral as
\begin{equation}
\label{eq:uniform-entropy-integral}
J(\delta)=\int_{0}^{\delta}\sup(1+\log N(\varepsilon\|F\|_{Q,2},\mathcal{F},L_{2}(Q)))^{1/2}d\varepsilon,
\end{equation}
where the supremum is taken over all probability measures $Q$ with $\|F\|_{Q,2}>0$.

We shall proceed in a way similar to that in \cite{CCK14-MR3262461}.
Let $M=\max\{F(X_{i}):1\leq i\leq n\}$.
Throughout this proof, $a\lesssim b$ denotes $a\leq Cb$ with $C$ depending only on $q$.
Observe that for $k\leq q$,
\begin{equation}
\label{eq:maximum-UB}
\|M\|_{k}\leq\|M\|_{q}\leq\Big(\sum_{i=1}^{n}\mathbb{E}|F(X_{i})|^{q}\Big)^{1/q}\leq n^{1/q}b.
\end{equation}
Note also that by Theorem 2.6.7 in \cite{VW23-MR4628026}, the covering number in \eqref{eq:uniform-entropy-integral} satisfies
\begin{equation}
\label{eq:covering-number-UB}
N(\varepsilon\|F\|_{Q,2},\mathcal{F},L_{2}(Q))\leq\frac{A}{3}\vc(\mathcal{F})(16e)^{\vc(\mathcal{F})}\Big(\frac{1}{\varepsilon^{1/2}}\Big)^{2\vc(\mathcal{F})}\leq\Big(\frac{A}{\varepsilon}\Big)^{v}.
\end{equation}
By integration by parts, we have for every $\delta\in(0,1)$,
\begin{equation}
\label{eq:uniform-entropy-integral-UB}
J(\delta)\leq2\delta\sqrt{v\log(A/\delta)},
\end{equation}
where $J(\delta)$ is defined in \eqref{eq:uniform-entropy-integral}.
Let $\varepsilon=\delta_{3}^{3/2}/n^{1/2}$.
Then
\begin{equation}
H_{n}(\varepsilon):=\log n\vee\log N(\varepsilon\|F\|_{P,2},\mathcal{F},L_{2}(P))\lesssim K.
\label{eq:pf-lem-CCK14-1}
\end{equation}

By the fact that $\sup_{f\in\mathcal{F}}P|f|^{3}\leq b^{3}\delta_{3}^{3}$ and Lemma 2.2 in \cite{CCK14-MR3262461}, it follows that
\begin{align}
\phi_{n}(\varepsilon)&:=\mathbb{E}\|\mathbb{G}_{n}\|_{\mathcal{F}_{\varepsilon}}\vee\mathbb{E}\|\mathbb{B}\|_{\mathcal{F}_{\varepsilon}}\lesssim bJ(\varepsilon)+\frac{bJ^{2}(\varepsilon)}{n^{1/2-1/q}\varepsilon^{2}}\notag\\&\lesssim\frac{b\delta_{3}^{3/2}K^{1/2}}{n^{1/2}}+\frac{bK}{n^{1/2-1/q}}\lesssim\frac{bK}{n^{1/2-1/q}},
\label{eq:pf-lem-CCK14-2}
\end{align}
that
\begin{align}
\mathbb{E}\|\mathbb{E}_{n}|f(X_{i})|^{3}\|_{\mathcal{F}}&\lesssim b^{3}\delta_{3}^{3}+\frac{b^{3/2}}{n^{1/2-3/2q}}\Big(b^{3/2}J(\delta_{3}^{3/2})+\frac{b^{3/2}J^{2}(\delta_{3}^{3/2})}{n^{1/2-3/2q}\delta_{3}^{3}}\Big)\notag\\&\lesssim b^{3}\delta_{3}^{3}+\frac{b^{3}\delta_{3}^{3/2}K^{1/2}}{n^{1/2-3/2q}}+\frac{b^{3}K}{n^{1-3/q}}\lesssim b^{3}\delta_{3}^{3}+\frac{b^{3}K}{n^{1-3/q}}=:\kappa^{3},
\label{eq:pf-lem-CCK14-3}
\end{align}
and that
\begin{equation}
\mathbb{E}\|\mathbb{G}_{n}\|_{\mathcal{F}\cdot\mathcal{F}}\lesssim b^{2}J(\delta_{4}^{2})+\frac{b^{2}J^{2}(\delta_{4}^{2})}{n^{1/2-2/q}\delta_{4}^{4}}\lesssim b^{2}\delta_{4}^{2}K^{1/2}+\frac{b^{2}K}{n^{1/2-2/q}}.
\label{eq:pf-lem-CCK14-4}
\end{equation}
By $\kappa\geq bK^{1/3}n^{1/q-1/3}$, for $q\in[4,\infty)$, $n\geq3$, and an absolute constant $c>0$,
\begin{align}
\delta_{n}(\varepsilon,\gamma)&=\frac{1}{4}P\Big(\Big(\frac{F}{\kappa}\Big)^{3}1\Big(\frac{F}{\kappa}>\frac{cn^{1/3}}{\gamma^{1/3}H_{n}^{1/3}(\varepsilon)}\Big)\Big)\lesssim\Big(\frac{b}{\kappa}\Big)^{q}\Big(\frac{\gamma^{1/3}H_{n}^{1/3}(\varepsilon)}{n^{1/3}}\Big)^{q-3}\notag\\&\lesssim\Big(\frac{n^{1/3-1/q}}{K^{1/3}}\Big)^{q}\Big(\frac{\gamma^{1/3}H_{n}^{1/3}(\varepsilon)}{n^{1/3}}\Big)^{q-3}\lesssim\frac{\gamma^{q/3-1}}{K}\lesssim\frac{\gamma^{q/3-1}}{\log n}\leq1.\label{eq:pf-lem-CCK14-5}
\end{align}
Note that $\delta_{n}(\varepsilon,\gamma)\equiv0$ for $q=\infty$, as $F/\kappa\leq b/\kappa\lesssim n^{1/3}/H_{n}^{1/3}(\varepsilon)$ in this case.

By Theorem 2.1 in \cite{CCK14-MR3262461},
\begin{equation}
\mathbb{P}(|Z-\tilde{Z}|\gtrsim\Delta_{n}(\varepsilon,\gamma))\lesssim\gamma(1+\delta_{n}(\varepsilon,\gamma))+\frac{\log n}{n},
\label{eq:pf-lem-CCK14-6}
\end{equation}
where, by combining \eqref{eq:pf-lem-CCK14-1}--\eqref{eq:pf-lem-CCK14-4},
\begin{align}
\Delta_{n}(\varepsilon,\gamma)&:=\phi_{n}(\varepsilon)+\frac{b\varepsilon}{\gamma^{1/q}}+\frac{b}{\gamma^{1/q}n^{1/2-1/q}}+\frac{b}{\gamma^{2/q}n^{1/2-1/q}}\notag\\&\quad+\frac{(\mathbb{E}\|\mathbb{G}_{n}\|_{\mathcal{F}\cdot\mathcal{F}})^{1/2}H_{n}^{1/2}(\varepsilon)}{\gamma^{1/2}n^{1/4}}+\frac{\kappa H_{n}^{2/3}(\varepsilon)}{\gamma^{1/3}n^{1/6}}\notag\\&\lesssim\frac{bK}{n^{1/2-1/q}}+\frac{b\delta_{3}^{3/2}}{\gamma^{1/q}n^{1/2}}+\frac{b}{\gamma^{1/q}n^{1/2-1/q}}+\frac{b}{\gamma^{2/q}n^{1/2-1/q}}\notag\\&\quad+\frac{b\delta_{4}K^{3/4}}{\gamma^{1/2}n^{1/4}}+\frac{bK}{\gamma^{1/2}n^{1/2-1/q}}+\frac{b\delta_{3}K^{2/3}}{\gamma^{1/3}n^{1/6}}+\frac{bK}{\gamma^{1/3}n^{1/2-1/q}}\notag\\&\lesssim\frac{bK}{\gamma^{1/2}n^{1/2-1/q}}+\frac{b\delta_{4}K^{3/4}}{\gamma^{1/2}n^{1/4}}+\frac{b\delta_{3}K^{2/3}}{\gamma^{1/3}n^{1/6}}=\Delta.
\label{eq:pf-lem-CCK14-7}
\end{align}
Then the desired result follows from combining \eqref{eq:pf-lem-CCK14-5}--\eqref{eq:pf-lem-CCK14-7}.

\section{Proof of \cref{thm:GA-mixing}}
\label{pf:thm-GA-mixing}
Unless explicitly stated, we shall denote $a\lesssim b$ denote $a\leq Cb$ for an absolute constant $C>0$ in the remaining proofs.
We first invoke the following coupling lemma that is a direct consequence of the \cite{Ber79-MR547109}'s lemma that can be found as Proposition 2 in \cite{DMR95-MR1324814}.

\begin{lemma}%[Implication of Berbee's lemma; ]
\label{lem:Berbee}
Let $(X_{i}:i\geq1)$ be a sequence of random variables taking values in a Polish space $\mathcal{X}$ defined on a probability space $(\Omega,\mathcal{A},\mathbb{P})$.
Suppose that $(\Omega,\mathcal{A},\mathbb{P})$ admits a sequence $(U_{i}:i\geq1)$ of uniform random variables on $(0, 1)$ independent of $(X_{i}:i\geq1)$.
Then there is a sequence $(X_{i}^{*}:i\geq1)$ of independent random variables such that $X_{i}^{*}=_{d}X_{i}$ and $\mathbb{P}(X_{i}^{*}\neq X_{i}) \leq \beta(\sigma(X_i), \sigma(X_j: j\geq i+1))$ for each $i\geq1$.
\end{lemma}

The following three lemmas are useful in the proof of \cref{thm:GA-mixing}.
In particular, \cref{lem:maximal-inequality-iid} is a maximal inequality for empirical processes indexed by classes of functions with i.i.d. observations, which is taken from Theorem 5.2 in \cite{CCK14-MR3262461}.
\cref{lem:maximal-inequality-naive} is a naive consequence of \cref{lem:maximal-inequality-iid} for a small fraction of empirical processes with general dependent observations.
\cref{lem:GA-blocking} is an extension of \cref{thm:CCK14} for empirical processes indexed by VC classes of functions with i.i.d. blocking observations.

\begin{lemma}
\label{lem:maximal-inequality-iid}
Let $X_{1},\dots,X_{n}$ be i.i.d. random variables taking values in a measurable space $(\mathcal{X},\mathcal{A})$ with a common distribution $P$ with the associated empirical process $\mathbb{G}_{n}=n^{-1/2}\sum_{i=1}^{n}(\delta_{X_{i}}-P)$.
Let $\mathcal{F}$ be a pointwise measurable class of measurable functions $\mathcal{X}\to\mathbb{R}$ with an envelope function $F$ satisfying $\|F\|_{P,2}<\infty$.
Suppose that there is a constant $\delta_{2}\in(0,1)$ such that $\sup\{\|f\|_{P,2}:f\in\mathcal{F}\}\leq\delta_{2}\|F\|_{P,2}$.
Then
\begin{equation}
\label{eq:lem-maximal-inequality-iid}
\mathbb{E}\|\mathbb{G}_{n}\|_{\mathcal{F}}\lesssim J(\delta_{2})\|F\|_{P,2}+\frac{J^{2}(\delta_{2})}{\delta_{2}^{2}n^{1/2}}\Big\|\max_{1\leq i\leq n}F(X_{i})\Big\|_{2}.
\end{equation}
\end{lemma}

\begin{lemma}%[Naive maximal inequality for a small fraction of empirical process indexed by classes of functions with general dependent observations]
\label{lem:maximal-inequality-naive}
Let the assumptions in \cref{lem:maximal-inequality-iid} be valid except that $X_{1},\dots,X_{n}$ are not necessarily independent observations (but still share a common distribution $P$).
Let $\mathbb{G}_{n,3}=n^{-1/2}\sum_{i=n-a_{3}+1}^{n}(\delta_{X_{i}}-P)$ be a small fraction of the summation of $\mathbb{G}_{n}$ for an integer $a_{3}\geq0$.
Then
\begin{equation}
\mathbb{E}\|\mathbb{G}_{n,3}\|_{\mathcal{F}}\lesssim\frac{a_{3}}{n^{1/2}}\Big(J(\delta_{2})+\frac{J^{2}(\delta_{2})}{\delta_{2}^{2}}\Big)\|F\|_{P,2}.
\end{equation}
\end{lemma}

\begin{lemma}%[Gaussian approximation for empirical processes indexed by VC classes of functions with i.i.d. blocking observations]
\label{lem:GA-blocking}
Let $Y_{1},\dots,Y_{m}$ be i.i.d. random variables taking values in a product measurable space $(\mathcal{X}^{a_{1}},\mathcal{A}^{a_{1}})$ with a common distribution $P_{1}$, $m\geq3$, and $a_{1}\geq1$.
Suppose that $Y_{1i}\sim P$ for each $1\leq i\leq a_{1}$.
Let $\mathcal{F}$ be a pointwise measurable class of measurable functions $\mathcal{X}\to\mathbb{R}$ with an envelope function $F$ such that $Pf=0$ for every $f\in\mathcal{F}$ and that $\vc(\mathcal{F})\leq v<\infty$.
Define a transformation $\varphi_{1}:\mathbb{R}^{\mathcal{X}}\to\mathbb{R}^{\mathcal{X}^{a_{1}}}$ with $f\mapsto(m/n)^{1/2}\sum_{i=1}^{a_{1}}f\circ\pi_{i}$ and $ma_{1}\leq n$ such that $\varphi_{1}(f)(Y_{1})=(m/n)^{1/2}\sum_{i=1}^{a_{1}}f(Y_{1i})$.
Suppose that there are constants $b>0$, $q\in[4,\infty]$, and $\delta_{3,1},\delta_{4,1}\in(0,1)$ such that $\|F\|_{P,q}\leq b$ and that $\sup\{\|\varphi_{1}f\|_{P_{1},k}:f\in\mathcal{F}\}\leq\delta_{k,1}\|\varphi_{1}F\|_{P_{1},k}$ for $k=3,4$.

Let $Z_{1}=\sup\{\mathbb{G}_{n,1}(\varphi_{1}f):f\in\mathcal{F}\}$, where $\mathbb{G}_{n,1}=m^{-1/2}\sum_{i=1}^{m}(\delta_{Y_{i}}-P_{1})$ is the empirical process associated to $Y_{1},\dots,Y_{m}$.
Let $G_{P_{1}}$ be a tight Gaussian process in $\ell^{\infty}(\varphi_{1}\mathcal{F})$ with mean zero and covariance function $\mathbb{E}(G_{P_{1}}(\varphi_{1}f)G_{P_{1}}(\varphi_{1}g))=P_{1}(\varphi_{1}f\varphi_{1}g)$ for every $f,g\in\mathcal{F}$.
Then for every $\gamma\in(0,1)$, there is a random variable $\tilde{Z}_{1}=_{d}\sup\{G_{P_{1}}(\varphi_{1}f):f\in\mathcal{F}\}$ such that
\begin{equation}
\mathbb{P}(|Z_{1}-\tilde{Z}_{1}|>\Delta_{1})\leq C_{q}\Big(\gamma+\frac{\log m}{m}\Big),
\end{equation}
where
\begin{equation}
\label{eq:Delta-1}
\Delta_{1}=\frac{bKn^{1/2}}{\gamma^{1/2}m^{1-1/q}}+\frac{b\delta_{4,1}K^{3/4}n^{1/2}}{\gamma^{1/2}m^{3/4}}+\frac{b\delta_{3,1}K^{2/3}n^{1/2}}{\gamma^{1/3}m^{2/3}},
\end{equation}
$K=c_{q}v\log(Am/\delta_{3,1}\delta_{4,1})$, $c_{q}$ and $C_{q}$ are positive constants depending only on $q$, and $A>0$ is an absolute constant.
\end{lemma}

\subsection{Proof of \cref{thm:GA-mixing}}
From the assumptions in \cref{thm:GA-mixing}, it can be shown that $\sup\{\|\varphi_{2}f\|_{P_{2},2}:f\in\mathcal{F}\}\leq\delta_{2}\|\varphi_{2}F\|_{P_{2},2}$, and $\sup\{\|\varphi_{1}f\|_{P_{1},k}:f\in\mathcal{F}\}\leq\delta_{k,1}\|\varphi_{1}F\|_{P_{1},k}$ for $k=3,4$, with $\delta_{k,1}=\delta_{k}$ for $k=3,4$.
Let $\varphi_{0}:\mathbb{R}^{\mathcal{X}}\to\mathbb{R}^{\mathcal{X}^{n}}$ with $f\mapsto n^{-1/2}\sum_{i=1}^{n}f\circ\pi_{i}$.
Rewrite $G_{P_{0}}$ as a tight Gaussian process in $\ell^{\infty}(\varphi_{0}\mathcal{F})$ with mean zero and covariance function
\begin{equation}
\sigma_{P_{0}}(f,g)=\mathbb{E}(G_{P_{0}}(\varphi_{0}f)G_{P_{0}}(\varphi_{0}g))=P_{0}(\varphi_{0}f\varphi_{0}g)=\mathbb{E}(\mathbb{G}_{n}f\mathbb{G}_{n}g).
\end{equation}
Similarly, let $G_{P_{1}}$ be a tight Gaussian process in $\ell^{\infty}(\varphi_{1}\mathcal{F})$ with mean zero and covariance function
\begin{equation}
\sigma_{P_{1}}(f,g)=\mathbb{E}(G_{P_{1}}(\varphi_{1}f)G_{P_{1}}(\varphi_{1}g))=P_{1}(\varphi_{1}f\varphi_{1}g)=\frac{ma_{1}}{n}\mathbb{E}(\mathbb{G}_{a_{1}}f\mathbb{G}_{a_{1}}g).
\end{equation}

In what follows, we shall use the big-block-small-block technique to decompose the proof into three steps.

\subsubsection*{Step 1: Approximating $\mathbb{G}_{n}$ by $\mathbb{G}_{n,1}(\varphi_{1}\cdot)$}
By \eqref{eq:Gn-decomposition}, it suffices to bound $\mathbb{G}_{n,2}(\varphi_{2}\cdot)$ and $\mathbb{G}_{n,3}$.
Recall that by \eqref{eq:uniform-entropy-integral-UB}, we have $J(\delta_{2})\lesssim\delta_{2}K^{1/2}$.
By \cref{lem:maximal-inequality-naive} and $n-(a_{1}+a_{2})m\leq2a_{1}$,
\begin{align}
\mathbb{E}\|\mathbb{G}_{n,3}\|_{\mathcal{F}}&\lesssim\frac{a_{1}}{n^{1/2}}\Big(J(\delta_{2})+\frac{J^{2}(\delta_{2})}{\delta_{2}^{2}}\Big)\|F\|_{P,2}\notag\\&\lesssim\frac{bn^{1/2}}{m}(\delta_{2}K^{1/2}+K)\lesssim\frac{bKn^{1/2}}{m}.
\end{align}
Then by Markov's inequality,
\begin{equation}
\label{eq:Gn-3-UB}
\mathbb{P}\Big(\|\mathbb{G}_{n,3}\|_{\mathcal{F}}\geq\frac{bKn^{1/2}}{\gamma m}\Big)\leq\gamma.
\end{equation}
Note that with a slight abuse of notation, the constants $K$, $c_{q}$, $C_{q}$, and $A$ in this proof may be different from case to case up to an absolute constant factor.

To bound $\mathbb{G}_{n,2}(\varphi_{2}\cdot)$, we invoke \cref{lem:Berbee} to reconstruct independent random variables $Y_{2}^{*},\dots,Y_{2m}^{*}$ such that $Y_{2i}^{*}=_{d}Y_{2i}$ and $\mathbb{P}(Y_{2i}^{*}\neq Y_{2i})\leq\beta(a_{1})$ for each $1\leq i\leq m$.
Let $\mathbb{G}_{n,2}^{*}=m^{-1/2}\sum_{i=1}^{m}(\delta_{Y_{2i}^{*}}-P_{2})$ be the empirical process associated to $Y_{2}^{*},\dots,Y_{2m}^{*}$.
Then we have
\begin{equation}
\label{eq:Gn-2-UB-1}
\mathbb{P}(\mathbb{G}_{n,2}\neq \mathbb{G}_{n,2}^{*})\leq\sum_{i=1}^{m}\mathbb{P}(Y_{2i}^{*}\neq Y_{2i})\leq m\beta(a_{1}).
\end{equation}
By \cref{lem:maximal-inequality-iid},
\begin{equation}
\label{eq:Gn-2-iid-UB-0}
\mathbb{E}\|\mathbb{G}_{n,2}^{*}\|_{\varphi_{2}\mathcal{F}}\lesssim J(\delta_{2})\|\varphi_{2}F\|_{P_{2},2}+\frac{J^{2}(\delta_{2})}{\delta_{2}^{2}m^{1/2}}\Big\|\max_{1\leq i\leq m}\varphi_{2}F(Y_{2i})\Big\|_{2}.
\end{equation}
For the first term in \eqref{eq:Gn-2-iid-UB-0}, by triangle inequality, we have
\begin{equation}
\label{eq:Gn-2-iid-UB-1}
\|\varphi_{2}F\|_{P_{2},2}\leq\frac{a_{2}m^{1/2}}{n^{1/2}}\|F\|_{P,2}\leq\frac{(a_{2}/a_{1})bn^{1/2}}{m^{1/2}}.
\end{equation}
For the second term, recall that by \eqref{eq:maximum-UB}, we have
\begin{equation}
\label{eq:Gn-2-iid-UB-2}
\Big\|\max_{1\leq i\leq m}\varphi_{2}F(Y_{2i})\Big\|_{2}\leq m^{1/q}\|\varphi_{2}F\|_{P_{2},q}\leq\frac{(a_{2}/a_{1})bn^{1/2}}{m^{1/2-1/q}}.
\end{equation}
Then by \eqref{eq:Gn-2-iid-UB-0}--\eqref{eq:Gn-2-iid-UB-2},
\begin{equation}
\mathbb{E}\|\mathbb{G}_{n,2}^{*}\|_{\varphi_{2}\mathcal{F}}\lesssim\frac{(a_{2}/a_{1})b\delta_{2}K^{1/2}n^{1/2}}{m^{1/2}}+\frac{(a_{2}/a_{1})bKn^{1/2}}{m^{1-1/q}}.
\end{equation}
Then by Markov's inequality,
\begin{equation}
\label{eq:Gn-2-UB-2}
\mathbb{P}\Big(\|\mathbb{G}_{n,2}^{*}\|_{\varphi_{2}\mathcal{F}}\geq\frac{b\delta_{2}K^{1/2}n^{1/2}}{(a_{1}/a_{2})\gamma m^{1/2}}+\frac{bKn^{1/2}}{\gamma m^{1-1/q}}\Big)\leq\gamma.
\end{equation}
Consequently, by \eqref{eq:Gn-3-UB}, \eqref{eq:Gn-2-UB-1}, and \eqref{eq:Gn-2-UB-2}, we have
\begin{align}
&\mathbb{P}\Big(\sup_{f\in\mathcal{F}}|\mathbb{G}_{n}f-\mathbb{G}_{n,1}(\varphi_{1}f)|\geq\frac{b\delta_{2}K^{1/2}n^{1/2}}{(a_{1}/a_{2})\gamma m^{1/2}}+\frac{bKn^{1/2}}{\gamma m^{1-1/q}}\Big)\notag\\&\leq\mathbb{P}(\mathbb{G}_{n,2}\neq \mathbb{G}_{n,2}^{*})+\mathbb{P}\Big(\|\mathbb{G}_{n,2}^{*}\|_{\varphi_{2}\mathcal{F}}\geq\frac{b\delta_{2}K^{1/2}n^{1/2}}{(a_{1}/a_{2})\gamma m^{1/2}}+\frac{bKn^{1/2}}{2\gamma m^{1-1/q}}\Big)\notag\\&\quad+\mathbb{P}\Big(\|\mathbb{G}_{n,3}\|_{\mathcal{F}}\geq\frac{bKn^{1/2}}{2\gamma m}\Big)\notag\\&\leq\gamma+m\beta(a_{1}).
\label{eq:pf-thm-GA-mixing-1}
\end{align}

\subsubsection*{Step 2: Approximating $\mathbb{G}_{n,1}$ by $G_{P_{1}}$}
By a similar argument as that in \eqref{eq:Gn-2-UB-1}, we can reconstruct independent random variables $Y_{1}^{*},\dots,Y_{2m-1}^{*}$ such that $Y_{2i-1}^{*}=_{d}Y_{2i-1}$ and $\mathbb{P}(Y_{2i-1}^{*}\neq Y_{2i-1})\leq\beta(a_{2})$ for each $1\leq i\leq m$.
Let $\mathbb{G}_{n,1}^{*}=m^{-1/2}\sum_{i=1}^{m}(\delta_{Y_{2i-1}^{*}}-P_{1})$ be the empirical process associated to $Y_{1}^{*},\dots,Y_{2m-1}^{*}$.
Then we have
\begin{equation}
\label{eq:Gn-1-coupling}
\mathbb{P}(\mathbb{G}_{n,1}\neq \mathbb{G}_{n,1}^{*})\leq\sum_{i=1}^{m}\mathbb{P}(Y_{2i-1}^{*}\neq Y_{2i-1})\leq m\beta(a_{2}).
\end{equation}

Let $Z_{1}=\sup\{\mathbb{G}_{n,1}(\varphi_{1}f):f\in\mathcal{F}\}$ and $Z_{1}^{*}=\sup\{\mathbb{G}_{n,1}^{*}(\varphi_{1}f):f\in\mathcal{F}\}$.
By \cref{lem:GA-blocking}, there is a random variable $\tilde{Z}_{1}^{*}=_{d}\sup\{G_{P_{1}}(\varphi_{1}f):f\in\mathcal{F}\}$ such that
\begin{equation}
\label{eq:Gn-1-GA}
\mathbb{P}(|Z_{1}^{*}-\tilde{Z}_{1}^{*}|>\Delta_{1})\leq C_{q}\Big(\gamma+\frac{\log m}{m}\Big),
\end{equation}
where $\Delta_{1}$ is defined in \eqref{eq:Delta-1}.
Then by \eqref{eq:Gn-1-coupling} and \eqref{eq:Gn-1-GA},
\begin{align}
\mathbb{P}(|Z_{1}-\tilde{Z}_{1}^{*}|>\Delta_{1})&\leq\mathbb{P}(\mathbb{G}_{n,1}\neq\mathbb{G}_{n,1}^{*})+\mathbb{P}(|Z_{1}^{*}-\tilde{Z}_{1}^{*}|>\Delta_{1})\notag\\&\leq C_{q}\Big(\gamma+\frac{\log m}{m}\Big)+m\beta(a_{2}).
\label{eq:pf-thm-GA-mixing-2}
\end{align}

\subsubsection*{Step 3: Approximating $G_{P_{1}}$ by $G_{P_{0}}$}
We shall invoke \cref{lem:comparison-process}.
By Dudley's inequality \citep[Corollary 2.2.9]{VW23-MR4628026}, we have for $i\in\{0,1\}$,
\begin{equation}
\mathbb{E}\|G_{P_{i}}\|_{(\varphi_{i}\mathcal{F})_{P_{i},\varepsilon}}\lesssim J(\varepsilon)\|\varphi_{i}F\|_{P_{i},2}.
\end{equation}
By a similar argument as that in \eqref{eq:Gn-2-iid-UB-1}, we have
\begin{equation}
\|\varphi_{0}F\|_{P_{0},2}\leq bn^{1/2}\quad\text{and}\quad\|\varphi_{1}F\|_{P_{1},2}\leq\frac{bn^{1/2}}{m^{1/2}}.
\end{equation}
Let $\varepsilon=m^{1/q-1}$.
Then by \eqref{eq:uniform-entropy-integral-UB}, we have $J(\varepsilon)\lesssim\varepsilon K^{1/2}$, and therefore
\begin{align}
b_{P_{0},P_{1}}(\gamma,\varepsilon)&=\sum_{i\in\{0,1\}}(\mathbb{E}\|G_{P_{i}}\|_{(\varphi_{i}\mathcal{F})_{P_{i},\varepsilon}}+\varepsilon\|\varphi_{i}F\|_{P_{i},2}(2\log2/\gamma)^{1/2})\notag\\&\lesssim\sum_{i\in\{0,1\}}(J(\varepsilon)+\varepsilon(2\log2/\gamma)^{1/2})\|\varphi_{i}F\|_{P_{i},2}\notag\\&\lesssim\Big(\varepsilon K^{1/2}+\frac{\varepsilon}{\gamma^{1/2}}\Big)bn^{1/2}\lesssim\frac{bK^{1/2}n^{1/2}}{\gamma^{1/2}m^{1-1/q}}.
\label{eq:pf-thm-GA-mixing-3-1-1}
\end{align}
Let $\beta=2\delta^{-1}\log(m\vee N_{P_{0},P_{1}}(\varepsilon))$.
Then it can be shown that $2\beta^{-1}\log N_{P_{0,}P_{1}}(\varepsilon)+3\delta\leq4\delta$ and $\eta\leq2m^{-1}\log m$.
Let
\begin{equation}
\delta=\frac{K^{1/2}}{\gamma^{1/2}}\|\sigma_{P_{0}}-\sigma_{P_{1}}\|_{\mathcal{F}^{2}}^{1/2}.
\label{eq:pf-thm-GA-mixing-3-1-2}
\end{equation}
Then by \eqref{eq:pf-thm-GA-mixing-3-1-1} and \eqref{eq:pf-thm-GA-mixing-3-1-2},
\begin{align}
\Delta_{P_{0},P_{1}}(\beta,\gamma,\delta,\varepsilon)&=2\beta^{-1}\log N_{P_{0,}P_{1}}(\varepsilon)+3\delta+b_{P_{0},P_{1}}(\gamma,\varepsilon)\notag\\&\lesssim\frac{K^{1/2}}{\gamma^{1/2}}\|\sigma_{P_{0}}-\sigma_{P_{1}}\|_{\mathcal{F}^{2}}^{1/2}+\frac{bK^{1/2}n^{1/2}}{\gamma^{1/2}m^{1-1/q}}.
\label{eq:pf-thm-GA-mixing-3-1}
\end{align}
By \eqref{eq:covering-number-UB}, we have $\log N_{P_{0},P_{1}}(\varepsilon)\lesssim K$ and therefore
\begin{equation}
\beta\delta^{-1}\|\sigma_{P_{0}}-\sigma_{P_{1}}\|_{\mathcal{F}^{2}}\lesssim K\delta^{-2}\|\sigma_{P_{0}}-\sigma_{P_{1}}\|_{\mathcal{F}^{2}}=\gamma.
\label{eq:pf-thm-GA-mixing-3-2}
\end{equation}
Consequently, by \eqref{eq:pf-thm-GA-mixing-3-1}, \eqref{eq:pf-thm-GA-mixing-3-2}, and \cref{lem:comparison-process}, there is a random variable $\tilde{Z}_{0}=_{d}\sup\{G_{P_{0}}(\varphi_{0}f):f\in\mathcal{F}\}$ such that
\begin{equation}
\mathbb{P}\Big(|\tilde{Z}_{1}^{*}-\tilde{Z}_{0}|>\frac{K^{1/2}}{\gamma^{1/2}}\|\sigma_{P_{0}}-\sigma_{P_{1}}\|_{\mathcal{F}^{2}}^{1/2}+\frac{bK^{1/2}n^{1/2}}{\gamma^{1/2}m^{1-1/q}}\Big)\lesssim\gamma+\frac{\log m}{m}.
\label{eq:pf-thm-GA-mixing-3}
\end{equation}
Finally, by combining \eqref{eq:pf-thm-GA-mixing-1}, \eqref{eq:pf-thm-GA-mixing-2}, and \eqref{eq:pf-thm-GA-mixing-3}, we complete the proof.

\subsection{Proof of \cref{lem:maximal-inequality-naive}}
Letting $n=1$ in \eqref{eq:lem-maximal-inequality-iid} of \cref{lem:maximal-inequality-iid}, we have
\begin{equation}
\mathbb{E}\|\delta_{X_{1}}-P\|_{\mathcal{F}}\lesssim\Big(J(\delta_{2})+\frac{J^{2}(\delta_{2})}{\delta_{2}^{2}}\Big)\|F\|_{P,2}.
\end{equation}
Then an application of triangle inequality leads to the desired result.

\subsection{Proof of \cref{lem:GA-blocking}}
Observe that
\begin{equation}
\|\varphi_{1}F\|_{P_{1},q}\leq a_{1}(m/n)^{1/2}\|F\|_{P,q}\leq(n/m)^{1/2}b.
\end{equation}
By \cref{thm:CCK14}, it suffices to check that $V(\varphi_{1}(\mathcal{F}))\leq v$.
We shall prove by contradiction.
Suppose that $V(\varphi_{1}(\mathcal{F}))>v$.
Then we can find $v+1$ points $(y_{1},t_{1}),\dots,(y_{v+1},t_{v+1})$ in $\mathcal{X}^{a_{1}}\times\mathbb{R}$ shattered by the subgraphs $\{(y,t):t<\varphi_{1}(f)(y)\}$ of $\varphi_{1}(\mathcal{F})$, namely for every subset $I\subset[v+1]$, there is a function $f\in\mathcal{F}$ with $t_{i}<(m/n)^{1/2}\sum_{j=1}^{a_{1}}f(y_{ij})$ for $i\in I$ and $t_{i}\geq(m/n)^{1/2}\sum_{j=1}^{a_{1}}f(y_{ij})$ for $i\notin I$.
This implies that $t_{i}^{\prime}:=(n/m)^{1/2}t_{i}-\sum_{j=2}^{a_{1}}f(y_{ij})<f(y_{i1})$ for $i\in I$ and $t_{i}^{\prime}\geq f(y_{i1})$ for $i\notin I$, namely $(y_{11},t_{1}^{\prime}),\dots,(y_{v+1,1},t_{v+1}^{\prime})$ in $\mathcal{X}\times\mathbb{R}$ are shattered by the subgraphs of $\mathcal{F}$.
Then we must have $\vc(\mathcal{F})\geq v+1$, which leads to a contradiction.
This completes the proof.

\section{Proof of \cref{lem:comparison-process}}
\label{pf:lem-comparison-process}
The proof of \cref{lem:comparison-process} is based on the following two lemmas.
\cref{lem:Strassen-thm} is an implication of Strassen's theorem that can be found as Lemma 4.1 in \cite{CCK14-MR3262461}.
\cref{lem:comparison-vector} is a coupling inequality for comparisons between Gaussian vectors.

\begin{lemma}%[Implication of Strassen's theorem; Lemma 4.1 in \citealp{CCK14-MR3262461}]
\label{lem:Strassen-thm}
Let $\mu$ and $\nu$ be Borel probability measures on $\mathbb{R}$, and let $V$ be a random variable defined on a probability space $(\Omega, \mathcal{A}, \mathbb{P})$ with distribution $\mu$.
Suppose that $(\Omega, \mathcal{A}, \mathbb{P})$ admits a uniform random variable on $(0,1)$ independent of $V$.
Let $\varepsilon>0$ and $\delta>0$ be two positive constants.
Then there is a random variable $W$ with distribution $v$ such that $\mathbb{P}(|V-W|>\delta) \leq \varepsilon$ if and only if $\mu(A) \leq \nu(A^\delta)+\varepsilon$ for every Borel subset $A$ of $\mathbb{R}$.
\end{lemma}

\begin{lemma}%[Coupling inequality for comparisons between Gaussian vectors]
\label{lem:comparison-vector}
Let $X=(X_{1},\dots,X_{n})\sim N(0,\Sigma_{X})$ and $Y=(Y_{1},\dots,Y_{n})\sim N(0,\Sigma_{Y})$.
Let $Z_{X}=\max_{1\leq i\leq n}X_{i}$.
Then for every $\beta>0$ and $\delta>1/\beta$, there is a random variable $\tilde{Z}_{Y}=_{d}\max_{1\leq i\leq n}\sum_{i=1}^{n}Y_{i}$ such that
\begin{equation}
\mathbb{P}(|Z_{X}-\tilde{Z}_{Y}|>2\beta^{-1}\log n+3\delta)\leq\frac{\varepsilon+C\beta\delta^{-1}\|\Sigma_{X}-\Sigma_{Y}\|_{\max}}{1-\varepsilon},
\end{equation}
where $\varepsilon=e^{-\alpha/2}(1+\alpha)^{1/2}<1$, $\alpha=\beta^{2}\delta^{2}-1>0$, and $C$ is an absolute constant.
\end{lemma}

Let $Z_{Q}=\sup\{G_{Q}(\varphi_{Q}f):f\in\mathcal{F}\}$, $N=N_{P,Q}(\varepsilon)$, $e_{\beta}=\beta^{-1}\log N$, and $b_{R}=\mathbb{E}\|G_{R}\|_{(\varphi_{R}\mathcal{F})_{R,\varepsilon}}+\varepsilon\|F_{R}\|_{R,2}(2\log2/\gamma)^{1/2}$, $R\in\{P,Q\}$.
Let $\{\varphi_{P}f_{1},\dots,\varphi_{P}f_{N_{P}}\}$ be an $\varepsilon\|F_{P}\|_{P,2}$-net of $(\varphi_{P}\mathcal{F},L_{2}(P))$ with $N_{P}=N(\varepsilon\|F_{P}\|_{P,2},\varphi_{P}\mathcal{F},L_{2}(P))$.
Let $\{\varphi_{Q}f_{N_{P}+1},\dots,\varphi_{Q}f_{N}\}$ be an $\varepsilon\|F_{Q}\|_{Q,2}$-net of $(\varphi_{Q}\mathcal{F},L_{2}(Q))$.
Then for each $R\in\{P,Q\}$ and every $f\in\mathcal{F}$, there is an $f_{i}\in\{f_{1},\dots,f_{N}\}$ such that $\|\varphi_{R}f-\varphi_{R}f_{i}\|_{R,2}<\varepsilon\|F_{R}\|_{R,2}$, or equivalently $\varphi_{R}f-\varphi_{R}f_{i}\in(\varphi_{R}\mathcal{F})_{R,\varepsilon}$.
By taking the $G_{R}$ operation and the supremum iteratively, $|Z_{R}-Z_{R}^{\varepsilon}|\leq\|G_{R}\|_{(\varphi_{R}\mathcal{F})_{R,\varepsilon}}$, where $Z_{R}^{\varepsilon}=\max_{1\leq i\leq N}G_{R}(\varphi_{R}f_{i})$.
By Borell's inequality \citep[Proposition A.2.1]{VW23-MR4628026}, we further have $\|G_{R}\|_{(\varphi_{R}\mathcal{F})_{R,\varepsilon}}\leq b_{R}$ with probability at least $1-\gamma/2$.
By Lemmas \ref{lem:Strassen-thm} and \ref{lem:comparison-vector} with $X=(G_{P}(\varphi_{P}f_{1}),\dots,G_{P}(\varphi_{P}f_{N}))\sim N(0,\Sigma_{X})$ and $Y=(G_{Q}(\varphi_{Q}f_{1}),\dots,G_{Q}(\varphi_{Q}f_{N}))\sim N(0,\Sigma_{Y})$, we have
\begin{equation}
\|\Sigma_{X}-\Sigma_{Y}\|_{\max}=\max_{1\leq i,j\leq N}|\sigma_{P}(f_{i},f_{j})-\sigma_{Q}(f_{i},f_{j})|\leq\|\sigma_{P}-\sigma_{Q}\|_{\mathcal{F}^{2}},
\end{equation}
and therefore for every Borel subset $A$ of $\mathbb{R}$,
\begin{align}
\mathbb{P}(Z_{P}\in A)&\leq\mathbb{P}(Z_{P}^{\varepsilon}\in A^{b_{P}})+\frac{\gamma}{2}\notag\\&\leq\mathbb{P}(Z_{Q}^{\varepsilon}\in A^{2e_{\beta}+3\delta+b_{P}})+\frac{\gamma}{2}+\frac{\eta+C\beta\delta^{-1}\|\sigma_{P}-\sigma_{Q}\|_{\mathcal{F}^{2}}}{1-\eta}\notag\\&\leq\mathbb{P}(Z_{Q}\in A^{2e_{\beta}+3\delta+b_{P}+b_{Q}})+\gamma+\frac{\eta+C\beta\delta^{-1}\|\sigma_{P}-\sigma_{Q}\|_{\mathcal{F}^{2}}}{1-\eta}.
\end{align}
By \cref{lem:Strassen-thm} again, we complete the proof.

\subsection{Proof of \cref{lem:comparison-vector}}
We proceed in a way similar to that in proving Theorem 4.1 of \cite{CCK14-MR3262461}.
Let $e_{\beta}=\beta^{-1}\log n$ and $Z_{Y}=\max_{1\leq i\leq n}\sum_{i=1}^{n}Y_{i}$.
By \cref{lem:Strassen-thm}, it suffices to show that for every Borel subset $A$ of $\mathbb{R}$,
\begin{equation}
\label{eq:pf-lem-comparison-vector-1}
\mathbb{P}(Z_{X}\in A)\leq\mathbb{P}(Z_{Y}\in A^{2e_{\beta}+3\delta})+\frac{\varepsilon+C\beta\delta^{-1}\|\Sigma_{X}-\Sigma_{Y}\|_{\max}}{1-\varepsilon}.
\end{equation}
For $x=(x_{1},\dots,x_{n})$, let $F_{\beta}(x)=\beta^{-1}\log(\sum_{i=1}^{n}e^{\beta x_{i}})$.
Then by (17) and Lemma 4.2 in \cite{CCK14-MR3262461}, there is a smooth function $g:\mathbb{R}\to\mathbb{R}$ such that $\|g^{\prime}\|_{\infty}\leq\delta^{-1}$, $\|g^{\prime\prime}\|_{\infty}\lesssim\beta\delta^{-1}$, and
\begin{align}
1_{A}\Big(\max_{1\leq i\leq n}x_{i}\Big)&\leq1_{A^{e_{\beta}}}(F_{\beta}(x))\leq\frac{1}{1-\varepsilon}g\circ F_{\beta}(x)\notag\\&\leq\frac{\varepsilon}{1-\varepsilon}+1_{A^{e_{\beta}+3\delta}}(F_{\beta}(x))\leq\frac{\varepsilon}{1-\varepsilon}+1_{A^{2e_{\beta}+3\delta}}\Big(\max_{1\leq i\leq n}x_{i}\Big).
\label{eq:pf-lem-comparison-vector-2}
\end{align}
Let $f=g\circ F_{\beta}$.
Then by \eqref{eq:pf-lem-comparison-vector-2}, we have
\begin{align}
\mathbb{P}(Z_{X}\in A)&\leq\frac{1}{1-\varepsilon}\mathbb{E}f(X)\leq\frac{1}{1-\varepsilon}\mathbb{E}f(Y)+\frac{1}{1-\varepsilon}|\mathbb{E}f(X)-\mathbb{E}f(Y)|\notag\\&\leq\mathbb{P}(Z_{Y}\in A^{2e_{\beta}+3\delta})+\frac{\varepsilon}{1-\varepsilon}+\frac{1}{1-\varepsilon}|\mathbb{E}f(X)-\mathbb{E}f(Y)|.
\label{eq:pf-lem-comparison-vector-3}
\end{align}
By \eqref{eq:pf-lem-comparison-vector-1} and \eqref{eq:pf-lem-comparison-vector-3}, it suffices to show that
\begin{equation}
\label{eq:pf-lem-comparison-vector-4}
|\mathbb{E}f(X)-\mathbb{E}f(Y)|\lesssim\frac{\beta}{\delta}\|\Sigma_{X}-\Sigma_{Y}\|_{\max}.
\end{equation}

To prove \eqref{eq:pf-lem-comparison-vector-4}, we define
\begin{equation}
h(x)=\int_{0}^{1}\frac{1}{2t}(\mathbb{E}f(\sqrt{t}x+\sqrt{1-t}Y)-\mathbb{E}f(Y))dt.
\end{equation}
Direct calculation gives
\begin{equation}
\partial_{i}h(x)=\int_{0}^{1}\frac{1}{2\sqrt{t}}\mathbb{E}\partial_{i}f(\sqrt{t}x+\sqrt{1-t}Y)dt,
\end{equation}
and
\begin{equation}
\partial_{i}\partial_{j}h(x)=\int_{0}^{1}\frac{1}{2}\mathbb{E}\partial_{i}\partial_{j}f(\sqrt{t}x+\sqrt{1-t}Y)dt.
\label{eq:pf-lem-comparison-vector-6}
\end{equation}
Note that $f$ is a smooth function.
By Lemma 1 in \cite{Mec09-MR2797946}, we have
\begin{equation}
\sum_{i=1}^{n}\mathbb{E}(X_{i}\partial_{i}h(X))=\sum_{i=1}^{n}\sum_{j=1}^{n}\Sigma_{X,ij}\mathbb{E}\partial_{i}\partial_{j}h(X),
\end{equation}
and
\begin{equation}
f(x)-\mathbb{E}f(Y)=\sum_{i=1}^{n}x_{i}\partial_{i}h(x)-\sum_{i=1}^{n}\sum_{j=1}^{n}\Sigma_{Y,ij}\partial_{i}\partial_{j}h(x).
\end{equation}
This implies that
\begin{equation}
\mathbb{E}f(X)-\mathbb{E}f(Y)=\sum_{i=1}^{n}\sum_{j=1}^{n}(\Sigma_{X,ij}-\Sigma_{Y,ij})\mathbb{E}\partial_{i}\partial_{j}h(X).
\end{equation}
By Lemma 4.3 in \cite{CCK14-MR3262461} and \eqref{eq:pf-lem-comparison-vector-6},
\begin{align}
|\mathbb{E}f(X)-\mathbb{E}f(Y)|&\leq\frac{1}{2}\|\Sigma_{X}-\Sigma_{Y}\|_{\max}\sum_{i=1}^{n}\sum_{j=1}^{n}\|\partial_{i}\partial_{j}f\|_{\infty}\notag\\&\leq\frac{1}{2}\|\Sigma_{X}-\Sigma_{Y}\|_{\max}(\|g^{\prime\prime}\|_{\infty}+2\|g^{\prime}\|_{\infty}\beta)\notag\\&\lesssim\frac{\beta}{\delta}\|\Sigma_{X}-\Sigma_{Y}\|_{\max}.
\end{align}
This completes the proof.

\section{Proofs of Lemmas \ref{lem:VC}--\ref{lem:cov-residual}}
\label{pf:lem-VC-cov-residual}
\subsection{Proof of \cref{lem:VC}}
\label{pf:lem-VC}
We only prove (i) as (ii) has been proven by \cref{fig:VC}.
Let $(x_{1},t_{1})$, $(x_{2},t_{2})$, $(x_{3},t_{3})$, and $(x_{4},t_{4})$ be four arbitrary points in $\mathcal{X}\times\mathbb{R}$.
We shall prove by contradiction.
Suppose that the subgraphs of $\mathcal{F}$ shatters these four points.
Suppose that $\mathcal{F}$ selects the single point $(x_{i},t_{i})$ at $s=s_{i}$ for $i\in[4]$, namely
\begin{equation}
\{(x_{i},t_{i})\}=\{(x,t):t<f(x,s_{i})\}\cap\{(x_{j},t_{j}):j\in[4]\}.
\end{equation}

We first claim that there is at most one point $(x_{i},t_{i})$ satisfying $F(x_{i})<s_{i}$.
Otherwise, suppose without loss of generality that $F(x_{1})<s_{1}$, $F(x_{2})<s_{2}$, and $F(x_{1})\leq F(x_{2})$.
Then $t_{1}<h(s_{1})$, $t_{2}<h(s_{2})$, $t_{1}\geq f(x_{1},s_{2})$, and $t_{2}\geq f(x_{2},s_{1})$.
As $F(x_{1})\leq F(x_{2})<s_{2}$, we must have $t_{1}\geq f(x_{1},s_{2})=h(s_{2})$.
Then $t_{2}<h(s_{2})\leq t_{1}<h(s_{1})\leq f(x_{2},s_{1})$.
Contradiction.

Next, as there must be at least three points $(x_{i},t_{i})$ satisfying $F(x_{i})\geq s_{i}$, we suppose without loss of generality that $F(x_{i})\geq s_{i}$ for $i\in[3]$ and $F(x_{1})\leq F(x_{2})\leq F(x_{3})$.
Then $t_{1}<g(x_{1},s_{1})$, $t_{2}<g(x_{2},s_{2})$, and $t_{3}<g(x_{3},s_{3})$.
As $F(x_{2})\geq F(x_{1})\geq s_{1}$, we must have $t_{2}\geq f(x_{2},s_{1})=g(x_{2},s_{1})$.
Then $g(x_{2},s_{1})\leq t_{2}<g(x_{2},s_{2})$.
This implies that $s_{1}<s_{2}$, and thus $t_{1}<g(x_{1},s_{1})\leq g(x_{1},s_{2})$.
However, $t_{1}\geq f(x_{1},s_{2})$.
Then we must have $s_{1}\leq F(x_{1})<s_{2}$.
Similarly, we have $s_{1}\leq F(x_{1})<s_{2}\leq F(x_{2})<s_{3}\leq F(x_{3})$ and $t_{3}\geq g(x_{3},s_{2})$.
Suppose that $\mathcal{F}$ selects the two-point set $\{(x_{1},t_{1}),(x_{3},t_{3})\}$ at $s=s_{13}$, namely
\begin{equation}
\{(x_{1},t_{1}),(x_{3},t_{3})\}=\{(x,t):t<f(x,s_{13})\}\cap\{(x_{i},t_{i}):i\in[4]\}.
\end{equation}

We then claim that $F(x_{1})<s_{13}$.
Otherwise, suppose that $F(x_{1})\geq s_{13}$.
Then $t_{2}\geq g(x_{2},s_{13})$ and $t_{3}<g(x_{3},s_{13})$.
By $g(x_{3},s_{2})\leq t_{3}<g(x_{3},s_{13})$, we must have $s_{2}<s_{13}$.
This implies that $t_{2}<g(x_{2},s_{2})\leq g(x_{2},s_{13})$.
Contradiction.

So far, we have not imposed any condition on $F(x_{4})$.
Next, we shall consider the two cases of $F(x_{4})<s_{4}$ and $F(x_{4})\geq s_{4}$.

Case 1: If $F(x_{4})<s_{4}$, then $t_{4}<h(s_{4})$.
Note also that $t_{1}\geq f(x_{1},s_{4})$ and $t_{4}\geq f(x_{4},s_{13})$.
However, by $F(x_{1})<s_{13}$, we must have $t_{4}<h(s_{4})\leq f(x_{1},s_{4})\leq t_{1}<f(x_{1},s_{13})=h(s_{13})\leq f(x_{4},s_{13})$.
Contradiction.

Case 2: If $F(x_{4})\geq s_{4}$, then we suppose that $\mathcal{F}$ selects the two-point set $\{(x_{2},t_{2}),(x_{4},t_{4})\}$ at $s=s_{24}$.
Suppose without loss of generality that $F(x_{1})\leq F(x_{2})\leq F(x_{3})\leq F(x_{4})$.
By the same argument, we must have $F(x_{2})<s_{24}$, $t_{2}<h(s_{24})$, and $t_{2}\geq f(x_{2},s_{13})$.
Then $F(x_{1})\leq F(x_{2})<s_{24}$.
This implies that $t_{1}\geq h(s_{24})$.
However, by $F(x_{1})<s_{13}$, we must have $t_{2}<h(s_{24})\leq t_{1}<h(s_{13})\leq f(x_{2},s_{13})$.
Contradiction.

To sum up, the subgraphs of $\mathcal{F}$ cannot shatter the aforementioned four points in either case.
This completes the proof.

\subsection{Proof of \cref{lem:cdf-rho}}
\label{pf:lem-cdf-rho}
The equation $\partial_{\rho}\Psi(x,y;\rho)=\psi(x,y;\rho)$ is a direct consequence of $\partial_{\rho}\Phi(x,y;\rho)=\phi(x,y;\rho)$.
An elegant derivation of $\partial_{\rho}\Phi(x,y;\rho)=\phi(x,y;\rho)$ can be found in (3)--(4) of \cite{Pla54-MR65047}.
This completes the proof.

\subsection{Proof of \cref{lem:cov-residual}}
\label{pf:lem-cov-residual}
By the decomposition that
\begin{equation}
\phi(x,y;\rho)=\frac{1}{\sqrt{1-\rho^{2}}}\phi(x)\phi\Big(\frac{y-\rho x}{\sqrt{1-\rho^{2}}}\Big),
\end{equation}
we have for $|x|,|y|\geq r$ and $|\rho|\leq\rho_{0}$,
\begin{align}
\frac{\phi(x,y;\rho)}{(\phi(x)\phi(y))^{1/2}}&=\frac{1}{\sqrt{1-\rho^{2}}}\frac{\phi^{1/2}(x)}{\phi^{1/2}(y)}\phi\Big(\frac{y-\rho x}{\sqrt{1-\rho^{2}}}\Big)\notag\\&=\frac{1}{(2\pi(1-\rho^{2}))^{1/2}}\exp\Big(\frac{y^{2}}{4}-\frac{x^{2}}{4}-\frac{(y-\rho x)^{2}}{2(1-\rho^{2})}\Big)\notag\\&=\frac{1}{(2\pi(1-\rho^{2}))^{1/2}}\exp\Big(-\frac{(1-\rho)^{2}(x^{2}+y^{2})+2\rho(x-y)^{2}}{4(1-\rho^{2})}\Big)\notag\\&\leq\frac{1}{(2\pi(1-\rho_{0}^{2}))^{1/2}}\exp\Big(-\frac{(1-\rho_{0})r^{2}}{2(1+\rho_{0})}\Big).\label{eq:pf-lem-cov-residual-1}
\end{align}
Note also that for $\pi_{0}(\lambda)\leq1/2$,
\begin{equation}
\sigma_{0}^{2}(\lambda)=\pi_{0}(\lambda)(1-\pi_{0}(\lambda))\geq1-\Phi(\lambda)\geq\frac{1}{\lambda+1/\lambda}\phi(\lambda)\geq\frac{1}{4\lambda}\phi(\lambda).
\label{eq:pf-lem-cov-residual-2}
\end{equation}
Then by \eqref{eq:pf-lem-cov-residual-1}--\eqref{eq:pf-lem-cov-residual-2}, for $\lambda,\nu\in[\lambda_{1},\lambda_{2}]$, we have
\begin{equation}
\frac{\phi(\lambda,\nu;\rho_{0})}{\sigma_{0}(\lambda)\sigma_{0}(\nu)}\leq\frac{4(\lambda\nu)^{1/2}\phi(\lambda,\nu;\rho_{0})}{(\phi(\lambda)\phi(\nu))^{1/2}}\leq\frac{2^{3/2}\lambda_{2}}{(\pi(1-\rho_{0}^{2}))^{1/2}}\exp\Big(-\frac{(1-\rho_{0})\lambda_{1}^{2}}{2(1+\rho_{0})}\Big).
\label{eq:pf-lem-cov-residual-3}
\end{equation}

Observe that by definition, we have
\begin{equation}
\frac{1}{(\log n)^{d}}=\pi_{0}(\lambda_{1})=2(1-\Phi(\lambda_{1}))\sim\frac{\sqrt{2}}{\sqrt{\pi}\lambda_{1}}e^{-\lambda_{1}^{2}/2},
\end{equation}
which implies that
\begin{equation}
\lambda_{1}=(2d\log\log n-\log\log\log n-\log(\pi d)+o(1))^{1/2}.
\label{eq:pf-lem-cov-residual-4}
\end{equation}
Similarly, we have
\begin{equation}
\lambda_{2}=(2\log n-(2c+1)\log\log n-\log\pi+o(1))^{1/2}.
\label{eq:pf-lem-cov-residual-5}
\end{equation}
Then by combining \eqref{eq:pf-lem-cov-residual-3} and \eqref{eq:pf-lem-cov-residual-4}--\eqref{eq:pf-lem-cov-residual-5}, we complete the proof.

\appendix
\section{A brief power analysis of the higher criticism}
\label{sec:DB}
For completeness, we provide an additional power analysis on the effect of $[\alpha_{1},\alpha_{2}]$ on the higher criticism test.
The analysis is mainly based on those in \cite{DJ04-MR2065195} and \cite{DHJ11-MR2815777}.
Let $\beta\in(1/2,1)$ and $r>0$.
For the setting in \cref{sec:intro} when $T_{j}$'s are independent $N(\mu_{j},1)$ variables, we further consider the local alternative
\begin{equation}
H_{1}:\quad\sum_{j=1}^{p}1(\mu_{j}\neq0)=p^{1-\beta}\quad\text{and}\quad|\mu_{j}|\in\{0,\pm(2r\log p)^{1/2}\}.
\end{equation}
Let $\HC_{p,\alpha}=(p\alpha(1-\alpha))^{-1/2}\sum_{j=1}^{p}(1_{\{p_{j}\leq\alpha\}}-\alpha)$ be the single-level higher criticism statistic at candidate significance level $\alpha\in(0,1)$.
When $\alpha=L_{p}/p^{s}$ for $s\in[0,1]$ and a polylogarithmic term $L_{p}$, it can be shown that
\begin{equation}
\rho_{s}(\beta)=\begin{cases}
(s^{1/2}-((1+s)/2-\beta)^{1/2})^{2}, & 1/2<\beta\leq(1+s)/2,\\
\infty, & (1+s)/2<\beta<1
\end{cases}
\end{equation}
is the detection boundary of $\HC_{p,\alpha}$ in the sense that
\begin{equation}
\mathbb{P}(\HC_{p,\alpha}>z_{\gamma})\to\begin{cases}
1, & r>\rho_{s}(\beta),\\
\gamma, & r<\rho_{s}(\beta),
\end{cases}
\end{equation}
where $z_{\gamma}$ is the upper $\gamma\in(0,1)$ quantile of $N(0,1)$.

It can be shown that the optimal detection boundary, namely the detection boundary of the higher criticism statistic $\HC_{p}^{*}$ over $\alpha\in[\alpha_{1},\alpha_{2}]=[1/p,1/2]$, is the lower envelope of the detection boundary of $\HC_{p,\alpha}$, given by
\begin{equation}
\rho^{*}(\beta)=\inf_{s\in[0,1]}\rho_{s}(\beta)=\begin{cases}
\beta-1/2, & 1/2<\beta<3/4,\\
(1-(1-\beta)^{1/2})^{2}, & 3/4<\beta<1.
\end{cases}
\label{eq:DB}
\end{equation}
When $[\alpha_{1},\alpha_{2}]=[L_{p}/p^{1-\eta},L_{p}/p^{1-\theta}]$, the detection boundary becomes
\begin{align}
&\rho_{\theta,\eta}^{*}(\beta)=\inf_{s\in[1-\theta,1-\eta]}\rho_{s}(\beta)\notag\\&=\begin{cases}
((1-\theta)^{1/2}-(1-\beta-\theta/2)^{1/2})^{2}, & 1/2<\beta\leq3/4-\theta/4,\\
\beta-1/2, & 3/4-\theta/4<\beta\leq3/4-\eta/4,\\
((1-\eta)^{1/2}-(1-\beta-\eta/2)^{1/2})^{2}, & 3/4-\eta/4<\beta\leq1-\eta/2,\\
\infty, & 1-\eta/2<\beta<1.
\end{cases}
\label{eq:DB-2}
\end{align}
Consequently, the exclusion of the range $\alpha\in[L_{p}/p^{1-\theta},1/2]$ makes the detection boundary \eqref{eq:DB-2} slightly deviates from the optimal one \eqref{eq:DB} for $\beta\in(1/2,3/4-\theta/4]$, indicating a minor power loss in this region and that the modified higher criticism test is still rate optimal.
However, the exclusion of the range $\alpha\in[1/p,L_{p}/p^{1-\eta}]$ can result in a severe power loss unless $\eta=0$.
In particular, \eqref{eq:DB-2} slightly deviates from \eqref{eq:DB} for $\beta\in(3/4-\eta/4,1-\eta/2]$, indicating that the modified test is still rate optimal in this region.
When signals are highly sparse with $\beta\in(1-\eta/2,1)$, we have $\rho_{\theta,\eta}^{*}(\beta)=\infty$, indicating that the modified higher criticism test is even not rate optimal in this region.
\bibliographystyle{imsart-nameyear} % Style BST file (imsart-number.bst or imsart-nameyear.bst)
\bibliography{MLT0122}       % Bibliography file (usually '*.bib')

%% or include bibliography directly:
% \begin{thebibliography}{}
% \bibitem{b1}
% \end{thebibliography}

\end{document}